\documentclass{article}
\usepackage[utf8]{inputenc}
\usepackage{cmap}

\usepackage{arxiv}

\usepackage[T1]{fontenc}    % use 8-bit T1 fonts
\usepackage{url}            % simple URL typesetting
\usepackage{booktabs}       % professional-quality tables
\usepackage[unicode]{hyperref}       % hyperlinks

\usepackage{termutf}
\usepackage{amsthm}
\usepackage{amssymb,amsmath,amsfonts,amscd}
\usepackage{graphicx}

\usepackage[matrix,arrow]{xy}
\usepackage[russian]{babel}

\usepackage{pstricks}

\usepackage{wasysym}
\usepackage{floatflt}
\usepackage{lipsum}

\title{Ideals of graphs: finding a set of generators}

\author{
  \fbox{Golod E.~S.}
   \And
 Osipov G.~A.\thanks{e-mail address: \texttt{georgy.osipov@graphics.cs.msu.ru}}
}

\begin{document}
\maketitle

\def\t{{\mathrm t}}
\def\d{{\mathrm d}}
\def\e{{\mathrm e}}
\def\x{{\mathbf x}}
\def\y{{\mathbf y}}
\def\t{{\mathbf t}}
\def\s{{\mathbf s}}
\def\r{{\mathbf r}}
\def\k{{\mathbf k}}
\def\u{{\mathbf u}}
\def\v{{\mathbf v}}
\def\m{{\mathbf m}}
\def\z{{\mathbf z}}
\def\d{{\partial}}
\def\mfm{{\mathfrak m}}
\def\mfH{{\mathfrak H}}

\begin{abstract}
In this paper, we consider homological properties of so-called graph ideals. Consider $\Gamma$ is a graph with vertices $t_1$, \ldots, $t_s$, without self-loops and multiple adjacencies. We can associate with such a graph an ideal $I(\Gamma)$ of polynomial ring
$$
A(\Gamma) = \Bbbk[t_1\ldots,t_s]
$$
over $\Bbbk$, generated by $x_{ij}=t_it_j$, $i\ne j$, corresponding to edges of $\Gamma$. 

The object of this paper is an algebra of Koszul homology
$$
H_\bullet((x_{ij}),A(\Gamma))
$$
of Koszul complex $K((\x_{ij}),A(\Gamma))$.

The result of this paper is finding a minimal multiplicative system of generators of this algebra for some graphs $\Gamma$.

There is an element $\starr$ in homology algebra corresponding each vertex in the graph, that should be included in every set of generators of each graph. This is a sufficient system for trees. Also, there is a generator element $\bigcirc$ for every cycle with length $n$ if $n\bmod3=2$. System of $\starr$-s and maybe $\bigcirc$ is sufficient for a graph with only one cycle. Also, here described a set of generators for a graph that is two cycles with exactly one common vertex. If a graph is two graphs with known algebra generators, connected by an edge, the answer for the whole graph is also described in this paper.

Bold letter $\x$ means a sequence $x_1$, \ldots. Same for other letters ($\t$, $\r$, $\u$, $\m$). This paper is done in Russian, further translation is expected.
\end{abstract}

% keywords can be removed

\newpage
\tableofcontents 

\newtheorem*{theor*}{Теорема}
\newtheorem{theor}{Теорема}
\newtheorem{exa}[theor]{Пример}
\newtheorem{lemm}[theor]{Лемма}
\newtheorem{thnote}[theor]{Замечание}
\newtheorem{cons}[theor]{Следствие}
\newtheorem{defin}[theor]{Определение}

\def\factor#1#2{\raisebox{3pt}{$#1$\kern-2pt}/\raisebox{-4.5pt}{\kern-2pt$#2$}}

\section{Предварительные сведения о комплексах и их гомологиях}
Этот параграф содержит в основном хорошо известные сведения. за деталями читатель отсылается, например, к \cite{dold}.

%\subsection{Комплексы и их гомологии}
Пусть $A$~--- кольцо, и для набора произвольных $A$-модулей $K_i$, $i\in \mathbb Z$ и гомоморфизмов 
$$
d_i:K_{i+1}\to K_{i}
$$
выполнено условие
$$
\Im d_{i+1} \subset \Ker d_{i},
$$
тогда такой набор модулей и гомоморфизмов $(K_\bullet,d_\bullet)$ называется комплексом. Отображение $d_i$ называется дифференциалом или граничным гомоморфизмом этого комплекса. Граничный гомоморфизм также может обозначаться $d(K)$, а его индекс допустимо опускать. Комплекс называют ограниченным, если лишь конечное число модулей $K_i$ отлично от нуля. Если $z\in K_i$, то $i$ называется размерностью $z$. Размерность $z$ будем обозначать $\dim z$.

Пусть дан комплекс $(K_\bullet,d_\bullet)$. Модуль $\Im d_{i-1}\subset K_i$ называется модулем границ и обозначается $B_i$. Модуль $\Ker d_{i}\subset K_i$ называется модулем циклов и обозначается $Z_i$. Поскольку имеется включение $B_i\subset Z_i$, имеет смысл говорить о фактормодуле $Z_i/B_i$, называемом $i$-тым модулем гомологий комплекса $K$. $i$-тый модуль гомологий обозначается $H_i$, а его элементы называются классами гомологий.

Морфизмом степени $s$ комплексов $(K_\bullet,d_\bullet)$ и $(C_\bullet,\d_\bullet)$ называется набор отображений 
$$
\alpha_i : K_i \to C_{i+s},
$$
такой, что для любого индекса $i$ имеет место коммутативная диаграмма
\begin{equation} \label{eq-compl-morph}
\begin{CD}
K_{i-1} @<d_i<< K_{i} \\
@VV\alpha_{i-1}V @VV\alpha_i V \\
C_{i+s-1} @<\d_{i+s}<< C_{i+s}.
\end{CD}
\end{equation}
Морфизм комплексов естественным образом индуцирует морфизм гомологий: $[\alpha_i]:H(K)_i\to H(C)_{i+s}$.
Если степень морфизма не указывается, будем считать её нулевой.

Также можно говорить о точной последовательности комплексов и морфизмов комплексов
$$
\ldots \to K_\bullet^{(i-1)} \to K_\bullet^{(i)} \to K_\bullet^{(i+1)} \to \ldots,
$$
что означает точность соответствующей последовательности модулей при любом $n$:
$$
\ldots \to K_n^{(i-1)} \to K_n^{(i)} \to K_n^{(i+1)} \to \ldots.
$$

Если дана короткая точная последовательность комплексов
$$
0\to K_\bullet' \xrightarrow{\alpha} K_\bullet \xrightarrow{\beta} K_\bullet'' \to 0,
$$
тогда имеет место следующая бесконечная точная последовательность гомологий этих комплексов, называемая длинной точной последовательностью гомологий:
$$
\ldots \to H''_{n+1} \xrightarrow{\gamma_{n+1}}
H'_n \xrightarrow{\overline\alpha_{n}} H_n \xrightarrow{\overline\beta_{n}} H''_n
\xrightarrow{\gamma_{n}} H'_{n-1} \to \ldots,
$$
в которой $\overline\alpha_\bullet$ и $\overline\beta_\bullet$~--- морфизмы, индуцированные, соответственно, $\alpha$ и $\beta$, а $\gamma_\bullet$~--- некоторый естественный гомоморфизм модулей, называемый связывающим. Напомним, как строится свзязывающий гомоморфизм
$$
H''_n \xrightarrow{\gamma_{n}} H'_{n-1}:
$$
для этого в схеме
$$
\begin{CD}
H_n'' @<<<
Z_n'' @>>>
K_n'' @<\beta_n<<
K_n    \\
@. @. @. @VV\d_n V \\
@. @. @.
K_{n-1} @<\alpha_{n-1}<<
K_{n-1}' @<<<
Z_{n-1}' @>>>
H_{n-1}'
\end{CD}
$$
берутся последовательно образы и прообразы.

\subsection{Тензорное произведение комлексов над полем}
Результаты, приведённые в данном разделе касаются тензорных произведений комплексов над полем.

Пусть даны комплексы $(K_\bullet,d_\bullet)$ и $(C_\bullet,\d_\bullet)$, состоящие из модулей над кольцом $A$. Их тензорным произведением над $A$ служит комлекс
$$
\ldots \leftarrow
(K\otimes_A C)_{n-1} \xleftarrow{(d\otimes\d)_{n}}
(K\otimes_A C)_{n}   \xleftarrow{(d\otimes\d)_{n+1}}
(K\otimes_A C)_{n+1} \leftarrow
\ldots,
$$
в котором $(K\otimes_A C)_i$ обозначает модуль
$$
\bigoplus_{k+l=i}(K_k\otimes_A C_l),
$$
а граничный гомоморфизм $(d\otimes\d)_i$ определяется на слагаемом $(K_k\otimes_A C_l)$ с условием $k+l=i$ следующим образом:
$$
(d\otimes\d)_i: z\otimes z' \mapsto d_k(z)\otimes z' + (-1)^{k} z\otimes \d_l(z').
$$

Гомологии тензорного произведения комплексов связаны с гомологиями этих комплексов посредством формулы Кюннета. Для дальнейших результатов понадобится частный случай этой формулы, верный, если тензорное произведение берётся над полем:
$$
H(K\otimes_\Bbbk C) = H(K)\otimes_\Bbbk H(C),
$$
где $(K_\bullet,d_\bullet)$ и $(C_\bullet,\d_\bullet)$~--- комплексы над некоторым полем $\Bbbk$.
В правой части этой формулы под $H(K)$ и $H(C)$ подразумеваются комплексы с нулевым дифференциалом, состоящие из модулей $H_i(K)$ и $H_i(C)$, а под $H(K)\otimes_\Bbbk H(C)$~--- тензорное произведение этих комплексов.

Учитывая, что над полем все комплексы свободны, с помощью формулы Кюннета можно получить изоморфизм 
$$
H(K\otimes C)\cong H(K)\otimes H(C),\;\factor{Z(K\otimes C)}{B(K\otimes C)}\cong H(K)\otimes H(C).
$$
Из сюрьективного отображения $Z(K)\otimes Z(C) \to H(K)\otimes H(C)$ получается тогда сюрьекция 
$Z(K)\otimes Z(C) \to \factor{Z(K\otimes C)}{B(K\otimes C)}$. Таким образом получается другая форма формулы Кюннета:
\begin{equation}\label{kunetform}
 Z(K\otimes_{\Bbbk}C)=Z(K)\otimes Z(C)+B(K\otimes C).
\end{equation}

Основная лемма, которая и используется при доказательстве теоремы о перешейке, формулируется следующим образом:
\begin{lemm}[о прообразе] \label{factorlemma}
Пусть $K$ и $C$ --- комплексы над полем, $K'\subseteq Z(K)$, $C'\subseteq Z(C)$ --- их подкомплексы с нулевым дифференциалом.

Тогда 
\begin{equation} \label{factorlemmaeq}
d^{-1}(K'\otimes C')=d^{-1}(K')\otimes C' + K'\otimes d^{-1}(C') + Z(K\otimes C).
\end{equation}
\end{lemm}

Ещё один интересный факт, который используется при доказательстве леммы о прообразе, сам по себе понадобится в дальнейшем, и поэтому выведен в отдельную лемму:

\begin{lemm}[о границах]\label{boardlemma}
Пусть $K$ и $C$ --- комплексы над полем, $K'\subseteq Z(K)$, $C'\subseteq Z(C)$ --- их подкомплексы с нулевым дифференциалом.

Тогда 
\begin{equation}\label{boardeq}
(K'\cap B(K))\otimes C' + K'\otimes(C'\cap B(C)) = (K'\otimes C')\cap B(K\otimes C)
\end{equation}
\end{lemm}
\prove{}{
Рассмотрим правую часть (\ref{boardeq}). Она является ядром отображения факторизации
$$
K'\otimes C'\xrightarrow{\varphi}H(K\otimes C)=H(K)\otimes H(C),
$$
которое является композицией проекции 
$$
K'\otimes C'\to \factor{K'}{B(K)\cap K'}\otimes \factor {C'}{B(C)\cap C'}
$$
и инъективного (как тензорное произведение инъективных) отображения 
$$
\factor{K'}{B(K)\cap K'}\otimes \factor {C'}{B(C)\cap C'} = 
\left(\factor{K'+B(K)}{B(K)}\right)\otimes\left(\factor{C'+B(C)}{B(C)}\right) \to 
H(K)\otimes H(C).
$$

Значит, ядро $\varphi$ равно 
\begin{equation}\label{boardker}
\Ker\left(K'\otimes C' \to \factor{K'}{B(K)\cap K'}\otimes \factor {C'}{B(C)\cap C'}\right).
\end{equation}

К (\ref{boardker}) применим известную формулу 
\begin{equation}
\factor M{M'}\otimes \factor N{N'}=\factor{M\otimes N}{M'\otimes N + M\otimes N'},
\end{equation}
верную для любых модулей $M$ и $N$ и их подмодулей $M'$ и $N'$,
и получим, что искомое ядро равно
\begin{multline*}
\Ker\left(K'\otimes C' \to \factor{K'\otimes C'}{(B(K)\cap K') \otimes C' + K'\otimes(B(C)\cap C')}\right)= \\
= (B(K)\cap K') \otimes C' + K'\otimes(B(C)\cap C'),
\end{multline*}
что и требовалось доказать.
}

Теперь осталось доказать лемму о прообразе.
\prove{ леммы \ref{factorlemma}}{
Напомним её формулировку:

$$
d^{-1}(K'\otimes C')=d^{-1}(K')\otimes C' + K'\otimes d^{-1}(C') + Z(K\otimes C). \eqno (\ref{factorlemmaeq})
$$

Очевидно, что 
$$
d^{-1}(K'\otimes C')=d^{-1}((K'\otimes C')\cap B(K \otimes C)),
$$
так что левую часть равенства (\ref{factorlemmaeq}) можно заменить на $d^{-1}(K'\otimes C'\cap B)$.
Поскольку $Z$ является ядром оператора $d$, чтобы установить справедливость (\ref{factorlemmaeq}) было бы достаточно доказать такой факт:
\begin{equation}\label{factorequiv}
d(d^{-1}(K')\otimes C' + K'\otimes d^{-1}(C'))=(K'\otimes C')\cap B(K \otimes C).
\end{equation}
Так как $d|_{K'}=0$, $d|_{C'}=0$, $d(d^{-1}(K'))=K'\cap B(K)$, а $d(d^{-1}(C'))=C'\cap B(C)$, левую часть равенства (\ref{factorequiv}) можно переписать в как
$$
(K'\cap B(K))\otimes C'+K'\otimes(C'\cap B(C)).
$$

Осталось только применить лемму \ref{boardlemma}.
}

%Докажем ещё один факт, который окажется полезным при изучении восьмёрок. 
%\begin{lemm}\label{caplemma}
%Если $M$ и $N$ --- векторные пространства над полем $\Bbbk$, $M'\subset M$, $N'\subset N$ подпространства, то 
%$$
%M'\otimes N\cap M \otimes N'=M'\otimes N'.
%$$
%\end{lemm}
%\prove{}{ Самое простое доказательство заключается в том, чтобы выбрать в $M$ и в $N$ по базису, согласованному, соответственно с $M'$ и $N'$, который будет обозначаться $m'_1,\ldots,m'_r$, $m_1,\ldots,m_s$ и $n'_1,\ldots,n'_p$, $n_1,\ldots,n_q$, с условиями $m'_i\in M'$, $m_i\in M-M'$, $n'_i\in N'$ и $n_i\in N-N'$. Тогда базис тензорного произведения есть 
%$$
%\{m'_i\otimes n'_j,m'_i\otimes n_j,m_i\otimes n'_j,m_i\otimes n_j\}.
%$$
%Ни один из элементов $M'\otimes N$ не содержит в своём разложении по базису элементов вида $m_i\otimes n'_j$ или $m_i\otimes n_j$, аналогичное утверждение есть про $N$. Поэтому если элемент содержится и в $M'\otimes N$ и в $M \otimes N'$, то он в разложении по базису содержит только элементы вида $m'_i\otimes n'_j$, а значит, лежит в $M'\otimes N'$.
%}

Пусть алгебра $R$ над кольцом $A$ представляется в виде прямой суммы своих $A$-подмодулей:
$$
R = R_0 \oplus R_1 \oplus \ldots \oplus R_n \oplus \ldots,
$$
причём если $x\in R_i$, а $y\in R_j$, то $xy\in R_{i+j}$. Тогда говорят, что на $R$ задана $\mathbb N$-градуировка, а $R$ называют градуированной алгеброй. Элементы подмодулей $R_i$ называются однородными.

Если $z\in R_i$ и $z\ne0$, то число $i$ будем называть размерностью $z$ и обозначать $\deg z$. Для градуированной алгебры выполнено следующее соотношение: если $zz'\ne0$, то $\deg zz'=\deg z+\deg z'$.

\subsection{Комплекс Козюля}
Градуированная $A$-алгебра $R$ называется дифференциальной, если на ней задан дифференциал $d:R\to R$, являющийся $A$-гомоморфизмом со свойствами:
\begin{itemize}
\item $d(R_n)\subset R_{n-1}$,
\item если $x\in R_i$, $y\in R_j$, то $d(xy) = d(x)y + (-1)^ixd(y)$.
\end{itemize}

Пусть $A$~--- комутативное кольцо с едицинцей и $\x=(x_1,\ldots,x_n)$~--- некоторый набор элементов $A$. Рассмотрим внешнюю алгебру свободного $A$-модуля $A^n$ с базисом $e_1$, \ldots, $e_n$. Определим на ней дифференциал 
$$
d(e_{i_1}\wedge\ldots\wedge e_{i_k}) = \sum\limits_j (-1)^{j+1} x_{i_j}\bigwedge\limits_{l\ne j}e_{i_l}.
$$
Такое определение $d$ превращает $\Lambda(A^n)$ в дифференциальную алгебру с дополнительным условием $d\circ d=0$, которое означает, что последовательность
$$
\Lambda_{n}(A^n)\to\Lambda_{n-1}(A^n)\to\ldots\to\Lambda_0(A^n)
$$
является ограниченным комплексом.

\begin{defin}[комплекса Козюля]
Полученный комплекс называется комплексом Козюля и обозначается 
$$
K(x_1,\ldots,x_n,A).
$$

Если дан $A$-модуль $M$, то вводится обозначение:
$$
K(x_1,\ldots,x_n,M)=K(x_1,\ldots,x_n,A)\otimes_A M
$$
($M$ в правой части этой формулы понимается как комплекс состоящий из одного ненулевого члена $M$ в нулевой компоненте).
\end{defin}

Операция умножения $\wedge$ в случае комплекса Козюля $K=K(x_1,\ldots,x_n,M)$ индуцирует операцию умножения классов гомологий, превращая градуированный модуль $H(K)$ в $A$-алгебру. Эта алгебра иначе обозначается $H(x_1,\ldots,x_n,M)$. Использоваться будут ещё обозначения $Z(x_1,\ldots,x_n,M)$ и $B(x_1,\ldots,x_n,M)$ для модулей циклов и границ комплекса $K$.%В силу того, что для элементов $z_1, z_2\in Z(K)$ выполнено $z_1\wedge z_2\in Z(K)$, а для элементов $z'_1,z'_2\in B(K)$ выполнено $z'_1\wedge$

Заметим, что по этому определению, если взят всего один элемент кольца $x\in A$, то комлекс $K(x,A)$ имеет вид
$$
0\leftarrow A \xleftarrow{x} A \leftarrow 0.
$$

Комплекс Козюля можно также рассматривать, как тензорное произведение над $A$, явным образом определённых, комплексов одного элемента:
$$
K(x_1,\ldots,x_n,A) = \bigotimes_i K(x_i,A).
$$

\subsection{Редукция}
\begin{lemm}[о редукции]\label{factor}
 Если $x_1,\ldots,x_m \in A$ регулярная последовательность, то существует изоморфизм алгебр
$$
  H(x_1,\ldots,x_n,A) \cong H(\ol{x_{m+1}},\ldots,\ol{x_n},A/(x_1,\ldots,x_m)).
$$
\end{lemm}
\prove{}%
{ 
  Индукция позволяет свести задачу к случаю $m=1$, а в этом случае утверждение регулярности означает, что $x_1$ не делитель нуля в $A$. Действительно, если утверждение доказано для некоторого $m$, то к началу регулярной последовательности $x_1,\ldots,x_m$, которое само является регулярной последовательностью можно применить утверждение индукции получив: $H(x_1,\ldots,x_n,A) \cong H(\ol{x_{m}},\ldots,\ol{x_n},A/(x_1,\ldots,x_{m-1}))$. Далее, воспользовавшись тем, что $x_m$ не делитель нуля в фактормодуле $A/(x_1,\ldots,x_{m-1})$ получаем требуемый изоморфизм $$H(\ol{x_{m}},\ldots,\ol{x_n},A/(x_1,\ldots,x_{m-1})) \cong H(\ol{x_{m+1}},\ldots,\ol{x_n},A/(x_1,\ldots,x_m)).$$

  Чтобы доказать базу, рассмотрим точную последовательность модулей
   $$
   0\to A \xrightarrow{x_1} A \to A/(x_1) \to 0.
   $$
   Она индуцирует точную последовательность комплексов
   $$
   0 \to K(\mathbf{x},A) \xrightarrow{x_1} K(\mathbf{x},A) \to  K(\mathbf{x},A/(x_1)) \to 0,
   $$
   а та, в свою очередь, точную последовательность модулей гомологий
   $$
   \cdots\to H_i(\mathbf{x},A) \xrightarrow{x_1} H_i(\mathbf{x},A) \to H_i(\mathbf{x},A/(x_1)) 
   \to H_{i-1}(\mathbf{x},A) \to \cdots.
   $$
   В которой левая стрелка, умножение на ${x_1}$ --- нулевое отображение. Следует это из того, что $H_i$ есть алгебра над кольцом $H_0=A/(\mathbf{x})$, которое аннулируется элементом $x_1$. Отсюда получаем короткую точную последовательность для каждого $i$:
   $$
     0 \to H_i(\mathbf{x},A) \to H_i(\mathbf{x},A/(x_1)) \to H_{i-1}(\mathbf{x},A) \to 0.
   $$
   Эта последовательность расщепляется, потому что к правому отображению можно предъявить обратное, которое работает по следующей схеме: данный цикл $z$ умножает на $e_1$. Таким образом, получаем цикл в $K(\mathbf{x},A/(x_1))$, т. к. $d(z\wedge e_1)=x_1z$. Можно проверить, что оно действительно будет обратным в нашем случае. Таким образом,
   $$
     H_i(\mathbf{x},A/(x_1))=H_i(\mathbf{x},A)\oplus H_{i-1}(\mathbf{x},A).
   $$
   Преобразуя левую часть получаем:
   \begin{multline*}
     H_i(\mathbf{x},A/(x_1))=H_i(0,x_2,\ldots,x_n,A/(x_1))=\\
     =(H(0,A/(x_1))\otimes H(x_2,\ldots,x_n,A/(x_1)))_i=\\
    =H_i(x_2,\ldots,x_n,A/(x_1))\oplus H_{i-1}(x_2,\ldots,x_n,A/(x_1)).
   \end{multline*}
   Из этого и из построения изоморфизма можно получить требуемое равенство.
}

Из определения $H(x_1,\ldots,x_n,M)$ сразу же получается следующее следствие:
\begin{cons}
В условиях предыдущей теоремы 
$$
H(x_1,\ldots,x_n,M) \cong H(\ol{x_{m+1}},\ldots,\ol{x_n},M/M(x_1,\ldots,x_m)).
$$
\end{cons}

\begin{lemm}[о вычислении образа при редукции]
Пусть класс гомологий $[z]\in H$ представлен в виде 
$$
[z] = [e_i\wedge z_1 + z_2],
$$
причём $z_2$ не содержит $e_i$. Тогда редукция $[z]$ по $e_i$ выражается следующей формулой:
\begin{equation}\label{reduction-fromula}
\ol{[z]} = \ol{[z_2]},
\end{equation}
в которой под $\ol{[z_2]}$ понимается образ $[z_2]$ при факторизации по $de_i$.
\end{lemm}
\prove{}
{ Следует из доказательства теоремы \ref{factor}.
}

\begin{defin}[обобщённой редукции]
Редукцию можно продолжить по формуле \ref{reduction-fromula} до морфизма комплексов $K(\Gamma)\to\ol{K(\Gamma)}$, где $\ol K$~--- редуцированный комплекс, то есть $K(\ol{x_{m+1}},\ldots,\ol{x_n},A/(x_1,\ldots,x_m))$. Образ $z\in K$ называется обобщённой редукцией элемента $z$.
\end{defin}

\begin{lemm}[пример доказательств гомологичности с помощью редукции]\label{reduction-example}
Класс $[z]$ равен нулю при выполнении любого из следующих условий:
\begin{itemize}
\item $z$ делится на $de_i$ для некоторого $i$,
\item $z$ делится на $e_i$ для некоторого $i$ (на самом деле, отсюда следует, что $z=0$),
\item $z$ представляется в виде 
$$
z = \sum\limits_j e_{i_j}\wedge\alpha_j
$$
для некоторых элементов $\alpha_j\in K$, если $e_{i_j}$ отвечают регулярному набору рёбер.
\end{itemize}
\end{lemm}
\prove{}
{ Во всех случаях, редукция $[z]$ по нужному регулярному набору обращается в нуль.
}

\begin{defin}[редукции класса гомологий, поднятия при редукции]
Редукцией класса гомологий $[z]\in H(K)$ называется его образ при изоморфизме, отвечающем редукции. Поднятием при редукции класса редуцированного комплекса называется такой класс исходного комплекса, редукция которого совпадает с данным.
\end{defin}

\begin{lemm}[о поднятии при редукции]\label{podniatie}
Пусть дан комплекс Козюля $K$ над некоторым кольцом $A$. Если обобщённая редукция по элементу $x_i\in A$ элемента комплекса $a\in K$ есть цикл $\overline{a}$, то тогда утверждается, что для некоторых $a_1$ и $a_2$ выполнено $da=e_i\wedge a_1 + x_i a_2$, и класс гомологий
$$
[a-e_i\wedge a_2]
$$ является поднятием $[\overline{a}]$ при такой редукции. $a_2$ в такой ситуации обозначается как $\frac{da}{x_i}$.
\end{lemm}
\prove{}
{ Так как $d\overline{a} = 0$  и обобщённая редукция~--- морфизм комплексов, обобщённая редукция $da$ равна 0. Равенство нулю $da$ в редуцированном комплексе влечёт представление его как $da=e_i\wedge a_1+x_i\cdot a_2$, где $a_2$ не содержит $e_i$. Если обозначить $a_2$ как $\frac{da}{x_i}$, то корректной является мнемоническая формула $\tilde a=\frac{d(e_i\wedge a)}{x_i}$, на самом деле означающая $\tilde a=a-e_i\wedge a_2$.

Очевидно, обощённая редукция этого элемента $\tilde a$ совпадает с обощённой редукцией $a$, значит, если бы $\tilde a$ был циклом, то это бы означало решение задачи. Рассмотрим $d\tilde a$: 
$$
d\tilde a=da-x_i\wedge a_2 + e_i\wedge da_2=e_i\wedge( a_1 + da_2 ).
$$

Равенство этого дифференциала нулю будет следовать из равенства нулю дифференциала элемента $e_i\wedge a_1+x_i\cdot a_2$, который является не просто циклом, но и границей:
$$
0=d(e_i\wedge a_1+x_i\cdot a_2)=e_i\wedge(-da_1)+x_i(a_1+da_2).
$$
Приравнивая к нулю коэффициент при $e_i$ получаем $da_1=0$, и следовательно, $x_i(a_1+da_2)=0$. А поскольку, $x_i$ не был делителем нуля, это влечёт $a_1+da_2=0$. 
}

\subsection{Точная последовательность Козюля, лемма о добавлении элемента} \label{addsection}
Для того, чтобы исследовать случай циклического графа нам понадобится ещё один инструмент, вытекающий из общей последовательности гомологий Козюля. Существование этой последовательности мы назовём леммой о добавлении элемента и докажем лишь в той общности, которая необходима для дальнейшего исследования.

\begin{lemm}[о добавлении элемента]\label{addlemma}
Пусть $\x$ обозначает последовательность $x_1$, \ldots, $x_n$ элементов некоторого кольца $A$, а $\x'$ есть $x_2$, $x_3$, \ldots, $x_n$. Тогда для любого $i$ имеет место точная последовательность гомологий Козюля:
\begin{equation}\label{exact}
H_i(\x',A) \xrightarrow{u} H_i(\x,A) \xrightarrow{v} (0:x_1)_{H_{i-1}(\x',A)} \to 0
\end{equation}
\end{lemm}
\prove{}{
Вспомним, что по определению комплекса Козюля имеет место равенство $K(\x,A)=K(x_1,A)\otimes K(\x',A)$. Комплекс $K(x_1,A)$ состоит из двух модулей, изоморфных $A$, поэтому $K_i(\x,A)$ можно записать и явно:
$$
K_i(\x,A)=K_i(\x',A)\oplus K_{i-1}(\x',A),
$$
причём дифференциал этого комплекса пару $(z_1,z_2)$ переводит в $(dz_1+(-1)^i x_1z_2,dz_2)$. Такая последовательность даёт точную последовательность модулей
$$
0 \to K_i(\x',A) \to K_i(\x,A) \to K_{i-1}(\x',A) \to 0,
$$
из которой можно получить гомологическую последовательность (не забывая о том, что дифференциал комплекса $K_{i-1}$ есть умножение на $x_1$):
$$
H_i(\x',A) \to H_i(\x,A) \to H_{i-1}(\x',A) \xrightarrow{x_1} H_{i-1}(\x',A).
$$

Если в этой последовательности заменить третий член на аннулятор умножения на $x_1$, получится требуемая последовательность (\ref{exact}).
}

Явный вид $u$ и $v$ можно получить анализируя доказательство леммы о добавлении элемента. Ответ состоит в том, что $u$ индуцировано простым вложением комплексов, а $v$ есть деление на $e_1$, и $v$ аннулирует элементы без $e_1$. Можно сформулировать этот результат, поскольку он понадобится в дальнейшем:

\begin{lemm}[явный вид отображений $u$ и $v$ при добавлении элемента]\label{uv-deletion-lemm}
1. Отображение $u$ из леммы о добавлении элемента переводит класс гомологий $[z']$, где $z'\in K(\x',A)$ в класс $[z]$, где $z\in K(\x,A)$ имеет такую же запись, как и $z'$.

2. Пусть $z\in Z(\x,A)$ цикл, 
\begin{equation}\label{e1z1plusz2}
z=e_1\wedge z_1+z_2,
\end{equation} 
где $z_1$ и $z_2$ не содержат $e_1$. Тогда $v([z])=[z_1]$\footnote{Здесь считаем, что $K(\x')\subset K(\x)$, и $z_1\in K(\x')$ так как $z_1$ не содержит $e_1$, так что $[z_1]\in H(\x')$.}.
\end{lemm}
\prove{}{
1. $u$ индуцированно вложением комплексов $K(\x',A)\to K(\x,A)$.

2. $v$ есть спуск отображения факторизации $K_i\to K_{i-1}(\x',A)$ на гомологии.
}

Так как мы условились считать, что $K(\x',A)\subset K(\x,A)$, если $z\in K(\x',A)$, запись $[z]$ можно понимать двояко, если рассматриваются оба комплекса. В случае, где непонятно, какой комплекс имеется ввиду, будем использовать запись $[z]_{H(\x')}$ для обозначаения образа $z$ при факторизации в $H(\x')$, и аналогично, запись $[z]_{H(\x)}$, так что $u([z]_{H(\x')}) = [z]_{H(\x)}$.% Если $[z]\in H(\x')$, то $u([z])$ будем обозначать $[z]_{H(\x)}$ или просто $[z]$ (поскольку эти элементы имеют одинаковую запись в разложении по базису).

\begin{lemm}[об образе произведения при добавлении элемента]
1. $u([z]\wedge[z'])=u([z])\wedge u([z'])$.
2. Пусть $z$ и $z'$ --- циклы из $K(\x,A)$. Тогда 
\begin{equation}\label{prodimage}
u\circ v([z]\wedge [z'])=u\circ v([z])\wedge[z']-(-1)^{\dim z}[z]\wedge u\circ v([z'])).
\end{equation}
\end{lemm}
\prove{}{
1. Из явного вида $u$.

2. Пусть $z$ и $z'$ раскладываются аналогично (\ref{e1z1plusz2}) следующим образом:
\begin{align*}
z &= e_1\wedge z_1+z_2;\\
z' &= e_1\wedge z'_1+z'_2,
\end{align*}
тогда $z\wedge z'=e_1\wedge (z_1\wedge z'_2+(-1)^{\dim z_2} z_2\wedge z'_1) + z_2\wedge z'_2$. Отсюда получаем: 
$$
v([z\wedge z'])=[z_1\wedge z'_2+(-1)^{\dim z_2} z_2\wedge z'_1]=[z_1\wedge z'_2+(-1)^{\dim z'_1} z'_1 \wedge z_2].
$$ 
Рассмотрим теперь правую часть (\ref{prodimage}). Её можно переписать следующим образом:
\begin{multline*}
[ e_1 \wedge ( (-1)^{\dim z_1} z_1 \wedge z'_1 - (-1)^{\dim z} z_1 \wedge z'_1 )
+ z_1\wedge z'_2+(-1)^{\dim z'_1} z'_1 \wedge z_2 ]=\\=[z_1\wedge z'_2+(-1)^{\dim z'_1} z'_1 \wedge z_2].
\end{multline*}
Осталось только применить $u$, чтобы левая и правая части оказались в одном комплексе.
}

\begin{cons}
Если $z'\in Z(\x')$, $[z]\in H(\x)$, то $v([z]\wedge [z'])=v([z]) \wedge [z']$\footnote{В данной формуле обозначение $[z']$ используется в двух смыслах: как обозначение элемента $H(\x)$ и элемента $H(\x')$}.
\end{cons}

\subsection{Мономиальный комплекс Козюля}
\begin{defin}[мономиального комплекса]
Если в качестве кольца $A$ выбрано кольцо многочленов над полем:
$$
A=\Bbbk[t_1,\ldots,t_s],
$$
а в качестве элементов для построения комплекса~--- приведённые мономы $x_1$, \ldots, $x_n$ из $A$, то комплекс Козюля
$$
K(x_1,\ldots,x_n,A)
$$
называется мономиальным.
\end{defin}

Помимо $\mathbb Z$-градуировки, которая будет называться гомологической и даёт структуру градуированной алгебры, в мономиальном комплексе вводится также более тонкая, мономиальная, градуировка. Мономиальная градуировка является $\Bbbk$-градуировкой каждого модуля $K_i$ как $\Bbbk$-векторного пространства, относительно неё дифференциал $d(K)$ является градуированным гомоморфизмом $K$.

Элемент 
$$
\mathfrak m=\alpha me_{i_1}\wedge \ldots \wedge e_{i_k}
$$ 
мономиального комплекса $K$ называется мономом, если $\alpha\in\Bbbk$, а $m\in A$~--- приведённый моном кольца многочленов. Его степенью называется моном кольца многочленов $\deg \mathfrak m = mx_{i_1}\cdot \ldots \cdot x_{i_k}$.

Элемент $z\in K\not=0$ называется однородным, если он является суммой мономов одинаковой степени. Эта степень (являющаяся мономом из $A$) называется тогда степенью $z$ и обозначается $\deg z$.

Если $z_1$ и $z_2$ однородные элементы мономиального комплекса $K$, и $z_1\wedge z_2\ne0$, то $\deg z_1\wedge z_2=\deg z_1\deg z_2$. Если $z$ однородный элемент $K$ и $dz\ne0$, то $\deg z=\deg dz$.

$K$ порождается своими мономами как векторное пространство. Поэтому любой элемент $z\in K$ можно представить как сумму однородных. Максимальные однородные слагаемые в таком представлении называются однородными компонентами $z$.

Если $z$ однородный элемент мономиального комплекса $K$ и $dz\ne0$, то моном $\deg z$ называется степенью класса гомологий $[z]$, и обозначается $\deg [z]$, а $[z]$ называется однородным классом гомологий.

Относительной степенью однородного элемента $z\in K$ по переменной $t_i$ (степенью $z$ относительно $t_i$) называется число $\deg_{t_i} z$, такое что $\deg z$ делится на $t_i^{\deg_{t_i} z}$, но не делится на $t_i^{1+\deg_{t_i} z}$.

Полной степенью $z$ называется число $|\deg z|$ равное сумме
$$
|\deg z| = \sum_i \deg_{t_i} z.
$$

Аналогичные определения даются для однородных классов гомологий.

\begin{exa}
Комплекс $K(xy,xz,\Bbbk[x,y,z])$ мономиальный. Элемент $\alpha=ze_{xy}-ye_{zx}$ однородный, его мультистепень равна $xyz$. $\deg_x \alpha=\deg_y \alpha=\deg_z \alpha=1$, $|\deg\alpha|=3$.
\end{exa}

\section{Комплекс Козюля графа}
\subsection{Определения и основные свойства комплекса Козюля графа}
\begin{defin}[комплекса Козюля графа]
Пусть дан граф $\Gamma$ без петель и кратных рёбер, $V(\Gamma)$ обозначает множество его вершин, а $E(\Gamma)$~--- множество его рёбер. Пусть $V(\Gamma)=\{t_1,\ldots,t_n\}$, и $\x = \{t_it_j|(t_i,t_j)\in E(\Gamma)\}$~--- набор бесквадратных мономов полной степени 2, отвечающих рёбрам графа. Обозначим $A=\Bbbk[\t]$. Тогда мономиальный комплекс Козюля
$$
K(\Gamma) = K(\x,A)
$$
называется комплексом Козюля графа $\Gamma$.
\end{defin}

\begin{exa}
Циклу из $n$ вершин отвечает комплекс Козюля $K(t_1t_2,t_2t_3,\ldots,t_{n-1}t_n,t_nt_1,\Bbbk[\t])$. Последовательности (то есть графу, получаемому из цикла удалением одного ребра) соответствует комплекс $K(t_1t_2,t_2t_3,\ldots,t_{n-1}t_n,\Bbbk[\t])$.
\end{exa}

\begin{lemm}[о гомологиях Козюля несвязного графа]
% Пусть граф $\Gamma$ разбивается на две несвязных компоненты $\Gamma_1$ и $\Gamma_2$, алгебры гомологий которых порождаются, соответственно, классами $[z_1]$, \ldots, $[z_k]$ и $[z'_1]$, \ldots, $[z'_l]$. Тогда образы классов $[\z]$ и $[\z']$ при гомоморфизмах, индуцируемых вложениями $\Gamma_1\to\Gamma$ и $\Gamma_2\to\Gamma$ порождают $H(\Gamma)$.
Если граф $\Gamma$ разбивается на несвязанные части $\Gamma_1$ и $\Gamma_2$ с комплексами, соответственно, $K_1$ и $K_2$, то 
$$
H=H_1\otimes_\Bbbk H_2.
$$
\end{lemm}
\prove{}
{ Обозначим вершины графа $\Gamma_1$ через $\t$, вершины графа $\Gamma_2$ через $\s$, а их рёбра, соответственно, через $\x$ и $\y$. Так как $\Bbbk[\t]=\Bbbk[\t,\s]/(\s)$, а $\Bbbk[\s]=\Bbbk[\t,\s]/(\t)$, выполнено
$$
\Bbbk[\s]\otimes_{\Bbbk[\t,\s]}\Bbbk[\t] = \Bbbk
$$
По определению комплекса Козюля
\begin{align*}
K(\x,\y,\Bbbk[\t,\s]) 
  & = K(\x,\Bbbk[\t,\s]) \otimes_{\Bbbk[\t,\s]} K(\y,\Bbbk[\t,\s]) \\
  & = K(\x,\Bbbk[\t])\otimes_{\Bbbk}\Bbbk[\s]\otimes_{\Bbbk[\t,\s]}\Bbbk[\t]\otimes_{\Bbbk}K(\y,\Bbbk[\s]) \\
  & = K(\x,\Bbbk[\t])\otimes_{\Bbbk}K(\y,\Bbbk[\s]).
\end{align*}

Теперь можно применить формулу Кюннета, которая в случае тензорного произведения над полем утверждает сущетсвование нужного изоморфизма.
}

\begin{cons}[о порождающих несвязного графа]
Пусть граф $\Gamma$ разбивается на несвязанные части $\Gamma$ и $\Gamma'$ с комплексами, соответственно, $K$ и $K'$. Если $H$ порождается набором $\z$, а $H'$~--- набором $\z'$, то $H(\Gamma)$ порождается объединением $\z$, $\z'$. При этом, если $\z$ и $\z'$ были минимальными системами порождающих, то их объединение есть минимальная система порождающих $H(\Gamma)$.
\end{cons}

\begin{defin}[неразложимого класса]
Однородный класс гомологий называется неразложимым, если он не содержится в подалгебре, порождённой классами меньшей абсолютной степени. 
\end{defin}

Вложению графов $\Gamma_1\to\Gamma_2$ отвечает соответсвующий морфизм их комплексов $K_1\to K_2$, которые поднимается до отображения гомологий $H_1\to H_2$ (не обязательно инъективного).

\begin{lemm}[о вложении полного подграфа]
Если $\Gamma_1$ полный подграф $\Gamma_2$, то вложению $\Gamma_1 \to \Gamma_2$ отвечает инъективное отображение гомологий $H_1\to H_2$.
\end{lemm}
\prove{}{
Пусть $z\in K(\Gamma_1)$. Поскольку подграф граф полный, это равносильно тому, что для любой вершины $v\in V(\Gamma_2)-V(\Gamma_1)$, $\deg_v z=0$. Если $[z]_{\Gamma_2}=0$, то $z = dz'$ для некоторого $z'\in K(\Gamma_2)$. Но тогда $\deg z=\deg z'$, а значит, $z'\in K(\Gamma_1)$, откуда следует $[z]_{\Gamma_1}=0$.
}

Обозначим через $\mfm(\Gamma)$ максимальный идеал $H(\Gamma)$, порождённой всеми элементами положительной абсолютной степени.

\begin{lemm}[о порождающих максимального идеала]\label{max-ideal-gen-lemm}
Пусть $([m_1],\ldots,[m_k])=\mfm$, тогда элементы 1, $[m_1]$, \ldots, $[m_k]$ порождают $H$ как алгебру.
\end{lemm}
\prove{}
{ Обозначим через $H'$ подалгебру, порождённую элементами $1$, $[m_1]$, \ldots, $[m_k]$ и предположим, что $H\not=H'$. Возьмём однородный класс $[z]\in H$ наименьшей абсолютной степени, не лежащий в $H'$. Так как $H-\mfm \subset H'$, он лежит в $\mfm$, поэтому выражается в виде
$$
[m] = \sum_i [\alpha_i][m_i],
$$
но $|\deg\alpha_i| < |\deg m|$, так что $[\alpha_i]\in H'$, что ведёт к противоречию.
}

Многие дальнейшие результаты будет намного удобнее формулировать пользуясь терминологией порождающих максимального идеала, а не алгебры.

\begin{defin}[несобственных классов; собственных классов; \#-класов]
Пусть дан граф $\Gamma$. Назовём \emph{несобственными классами гомологий} классы из $H(\Gamma)$, лежащие в идеале, порождённом всевозможыми образами $v(\mfm(\Gamma'))$, где $\Gamma'$~--- некоторый собственный подграф $\Gamma$, а $v$~--- вложение гомологий, индуцированное вложением графов $\Gamma'\to\Gamma$. Идеал неcобственных классов обозначим $\mfH$. Элементы факторалгебры $H/\mfH$, обозначаемой $H_\#$ называются \SH-классами. Образ $z\in Z$ в $H/\mfH$ обозначается $[z]^\#$. Собственными классами называются такие класссы $[z]$, что $[z]^\#\ne0$
\end{defin}

В большинстве случаев, при изучении какого-либо класса графов все его подграфы уже хорошо изучены, поэтому в первую очередь, интерес будут представлять \SH-классы. 

\begin{lemm}[о степени собственных классов графа]\label{sobstven-lemm}
Пусть $[z]\in H$~--- однородный класс гомологий мультистепени $m$ графа $\Gamma$. Если $[z]^\#$~--- собственный, то $m$ содержит все переменные в положительной степени.
\end{lemm}
\prove{}
{ Пусть $m$ не содержит какую-либо переменную $t_i$. Рассмотрим граф $\Gamma'=\Gamma-t_i$. Тогда можно рассмотреть элемент $z_{\Gamma'}\in K'$, записывающийся также, как и $z$. Образ $[z_{\Gamma'}]$ при вложении гомологий $\Gamma'\to\Gamma$ равен $[z]$, значит, $[z]\in\mfH$.
}

\begin{lemm}[о сведении к факторалгебре]
Пусть классы гомологий $[z_1]$, \ldots, $[z_k]$ порождают $\mfH$, а образы классов $[z_{k+1}]$, \ldots, $[z_l]$ порождают максимальный идеал факторалгебры $H/\mfH$. Тогда $[z_1]$, \ldots, $[z_l]$ порождают $\mfH$.
\end{lemm}
\prove{}
{ Максимальный идеал $H/\mfH$ совпадает с $\mfm/\mfH$, поэтому элементы $[\z]$ обладают тем свойством, что порождают $\mfH$, а их образы порождают $\mfm/\mfH$.
}

Чтобы применять редукцию в комплексе графа, требуется определять, какие последовательности являются регулярными. Поэтому без доказательства (в силу своей простоты) приведём следующий факт:
\begin{lemm}[о регулярном наборе рёбер]
В комплексе графа регулярной последовательностью элементов является набор рёбер, никакие два из которых не имеют общей вершины.
\end{lemm}

\subsection{Звёздочки. Идеал вершины}
\begin{defin}
\emph{Звёздочкой на рёбрах $xy_1$, $xy_2$, \ldots, $xy_k$}, где $x$ и $y_i$ --- вершины графа, будет называться следующий класс гомологий:
$$
\starr_x(y_1,\ldots,y_k)=\left[\dfrac{d(e_{xy_1}\wedge\ldots\wedge e_{xy_k})}{x}\right],
$$
где деление на $x$ означает сокращение всех коэффициентов на $x$. Звёздочкой также будет называться граф состоящий из набора рёбер с общей вершиной, то есть такой, для которого звёздочка является собственным классом.
\end{defin}

При изменении порядка рёбер $xy_1$, $xy_2$, \ldots, $xy_k$, звёздочка меняет знак в случае, если перестановка была нечётной. При чётной перестановке рёбер, звёздочка не меняется.

Приведём явный вид звёздочки:
$$
\starr_x(y_1,\ldots,y_k)=\left[-\sum\limits_i (-1)^{i}y_i\bigwedge\limits_{j\not=i}e_{xy_j}\right].
$$

\begin{lemm} 
Звёздочки является нетривиальным неразложимым классом гомологий.
\end{lemm}
\prove{}
{ Докажем вначале, что не может быть цикла с постоянными коэффициентами (то есть такого, что $|\deg z|=2\dim z$). Действительно, его однородными компонентами были бы мономы, которые в исходном комплексе не являются циклами. Так что коэффициенты цикла --- многочлены, степени по крайней мере 1.

Равенство $d(\starr)=0$ следует из того, что $dd=0$. Все коэффициенты звёздочки --- многочлены степени 1, но у границ и у произведений элементов $Z$ коэффициенты имеют степень по крайней мере 2, так что звёздочка не может порождаться границами и разложимыми элементами, а следовательно, она нетривиальна и неразложима.
}

\begin{exa}
Звёздочка на одном ребре есть переменная: $\starr_x(y)=[y]$. Звёздочка на двух рёбрах $xa$ и $xb$ есть соотношение между $e_{xa}$ и $e_{xb}$:
$$
\starr_x(a,b)=[be_{xa}-ae_{xb}].
$$
\end{exa}

%Класс $[t_i]$ считаем звёздочкой даже в том случае, когда $t_i$ соответствует вершине степени 0.

Выпишем также относительные степени звёздочки, относительно вершин, которые в ней участвуют:
\begin{align*}
\deg_x    \starr_x(y_1,\ldots,y_k) & =k-1, \\
\deg_{y_i}\starr_x(y_1,\ldots,y_k) & =1.
\end{align*}

\begin{lemm}[о редукции звёздочек]
Пусть звёздочка $\starr_y(x_1,\ldots,x_s)$ редуцируется по ребру $xy_i$. Тогда её образ равен
$$
\ol{\starr_x(y_1,\ldots,y_s)} = [y_i\bigwedge\limits_{j\not=i}e_{xy_j}].
$$
\end{lemm}

\begin{defin}[идеала вершины]
Идеалом вершины $t_i$ графа $\Gamma$ с комплексом $K(\Gamma)$ называется идеал алгебры гомологий этого комплекса $H(\Gamma)$, порождённый всеми звёздочками в вершине $t_i$. Обозначение: $(\starr_{t_i})$. Идеал, порождённый всеми звёздочками в вершине $t_i$, содержащими рёбра $t_is_1$, \ldots, $t_is_l$ будет обозначаться $\starr_{t_i}(\s,\ldots)$.
\end{defin}
Далее будет дана другая характеризация идеала вершины

%Ближайшие результаты покажут, что совокупность всех звёздочек иногда является системой, порождающей всю алгебру гомологий графа. Так происходит в большинстве случаев. Однако попытки доказать этот факт в общем случае оканчиваются неудачей, что даёт повод ввести следующие понятия:

\begin{defin}[идеала звёздочек; \S-разложимых элементов; \SS-классов]
Сумма всех идеалов вершин называется идеалом звёздочек. Для него будет принято следующее обозначение:
$$
(\starr)=\sum\limits_i (\starr_{t_i}).
$$
Класс гомологий, лежащий в идеале звёздочек называется \S-разложимым. Класс, не являющийся \S-разложимым будет называться \S-неразложимым. Элементы факторалгебры $H/(\starr)$ называются \SS-классами (незвёздочными классами).
\end{defin}

%Далее (а именно, в параграфе \ref{ann-sect}) будет дана более простая характеризация идеала вершины. Также будет показано (см. \ref{kruzochki-sect}), что $\starr$-неразложимый класс не обязан быть неразложимым. 

%Можно также рассматривать факторалгебру $H$ по идеалу звёздочек $(\starr)$, которая обозначается $H_*=H/(\starr)$. Образ класса гомологий $[z]$ при такой факторизации будет обозначен $[z]_*$ или $z_*$, и назван $\starr$-классом цикла $z$.

\begin{lemm}[о линейных зависимостях среди звёздочек]
Пусть звёздочка $\starr_0$ представлена в виде минимальной однородной линейной комбинации с ненулевыми коэффициентами из поля $\Bbbk$ других звёздочек, $\starr_1$, \ldots, $\starr_k$, не равных $\starr_0$ и $-\starr_0$. Тогда $\starr_0$ есть звёздочка на двух рёбрах $\alpha\beta_1$ и $\alpha\beta_2$, вершины $\beta_1$ и $\beta_2$ соединены ребром, $k=2$, и выполнено:
\begin{align*}
\starr_1 &= \pm\starr_{\beta_1}(\alpha,\beta_2),\\
\starr_2 &= \pm\starr_{\beta_2}(\alpha,\beta_1),
\end{align*}
либо наоборот. В указанных условиях линейная зависимость существует.
\end{lemm}
\prove{}
{ Пусть звёздочка $\starr_0$ является звёздочкой в вершине $\alpha$ более, чем на двух рёбрах. В этом случае вершина $\alpha$ однозначно определяется по мультистепени соотношением $\deg_\alpha\starr_0>1$. Остальные вершины звёздочки также могут быть однозначно установлены, поскольку для них выполнено равенство $\deg_\beta\starr_0=1$. Поэтому получаем, что в данной мультистепени имеется только одна звёздочка (с точностью до знака). Совсем прост случай звёздочки на одном ребре, ведь наличие тогда линейной зависимости, означало бы линейную зависимость между переменными, которой не существует.

Остался только один случай, когда $\starr_0$ --- звёздочка на двух рёбрах. Обозначим эти рёбра как $\alpha\beta_1$ и $\alpha\beta_2$. Тогда $\deg\starr_0=\alpha\beta_1\beta_2$. Другая линейно-независимая звёздочка, очевидно, существует только в случае наличия ребра между $\beta_1$ и $\beta_2$. Таких звёздочек будет две (одна в вершине $\beta_1$, другая --- в $\beta_2$), и они описаны в условии леммы. Редукция по ребру $\beta_1\beta_2$ позволяет сразу же установить заявленную линейную зависимость.
}

Доказанная лемма позволяет находить \textit{минимальную систему звёздочек}, то есть минимальный набор звёздочек, порождающий идеал $(\starr)$. Для этого из множества всех звёздочек нужно исключить по одной звёздочке для каждого цикла длины три в графе.

\subsection{Деревья. Лемма об ограничении степени}
\begin{lemm}[об ограничении относительной степени]\label{ogranichenie-theor}
Пусть $v$ --- вершина графа $\Gamma$, и $\deg_{\Gamma} v$ --- валентность этой вершины (количество инцидентных с ней рёбер), $[z]\in H$~--- собственный однородный класс гомологий графа~$\Gamma$. Пусть также $u_1$, \ldots, $u_{\deg_{\Gamma} v}$ — вершины, с которыми соединена $v$. Тогда
\begin{enumerate}
\item Если $\deg_v z > \deg_{\Gamma} v$, то $[z]$ делится на $v$, и следовательно разложим.
\item Если $\deg_v z = \deg_{\Gamma} v$, то 
$$
[z]\in \bigcap\limits_{1\leqslant i\leqslant\deg_{\Gamma} v}(\starr_{u_i}(v,\ldots)),
$$ %представляется в виде
%$$
%[z]=\sum_i [z_i]\wedge\starr_{u_i}(v,v_{1,i},\ldots,v_{k,i})
%$$
%для некоторых $z_i\in Z$ и вершин $u_i$, $v_{j,i}$. То есть является линейной комбинацией разложимых элементов и, возможно, звёздочек.
\end{enumerate}
то есть является $\starr$-разложимым.
\end{lemm}
\newtheorem*{hypoteza*}{Гипотеза}
\begin{hypoteza*}
Если в предположениях леммы \ref{ogranichenie-theor} выполнено $\deg_v z = \deg_{\Gamma} v$, то верно включение 
$$
[z]\in \prod\limits_{1\leqslant i\leqslant\deg_{\Gamma} v}(\starr_{u_i}).
$$
\end{hypoteza*}

\prove{ леммы об ограничении степени}
{ \par\indent 1. Пусть $\deg_v[z]>\deg v$. Докажем, что $[z]=v[z_1]$, для некоторого $z_1\in Z$, т. е. $[z]$ разложим. Для этого достаточно показать, что любой моном $z$ делится на $v$. Пусть $m=\alpha\wprod i1k$ --- моном $z$. 
$$\deg v<\deg_vz=\deg_v m=\deg_v\alpha+\deg_v\wprod i1k.$$
Теперь, чтобы получить требуемое неравенство, $\deg_v\alpha>0$, осталось доказать, что $\deg_v\wprod i1k\leqslant\deg v$, но это прямо следует из того, что $\deg_v e_i\leqslant 1$, причём равенство выполнено только для $e_i$, отвечающих рёбрам, инцидентным с вершиной $v$, а таких рёбер всего $\deg v$. В данном случае $0=dz=d(vz')=vdz'$ влечёт $dz'=0$, т. е. $z'\in Z$, что и требуется.

\newbox\grpic\setbox\grpic=\hbox{\includegraphics{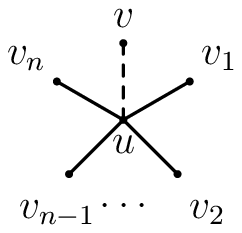}}
\begin{floatingfigure}[l]{\wd\grpic}\copy\grpic\end{floatingfigure}
\par\indent 2. %Теперь настало время применить редукцию. Редукция чрезвычайно полезна для исследования звёздочек, поскольку редукция по ребру, на котором строится звёздочка превращает её в моном. 
Применим редукцию, чтобы превратить звёздочку в класс монома. Кроме того, с помощью редукции можно понизить степень вершины, и это позволит применить рассуждение, аналогичное доказательству пункта 1, однако в неизменном виде оно не проходит, потому что в редуцированном комплексе существуют делители нуля, и это обстоятельство не позволяет вывести $dz'=0$ из утверждения $vdz'=0$. 

Итак, пусть $vu$ --- любое ребро, выходящее из вершины $v$. Докажем, что $[z]\in (\starr_u)$, то есть 
$$
[z]=\sum_i [z_i]\wedge\starr_{u}(v,v_{1,i},\ldots,v_{k,i})
$$
для некоторых циклов $z_i\in Z$ и вершин $v_{j,i}$.

Производя редукцию по ребру $uv$, получаем комплекс $K'$ над кольцом $A'=A/(vu)$, с условиями $H'=H$, $vu=0$ и $\deg'v=\deg v-1$. Поскольку теперь имеет место неравенство $\deg' v<\deg_v z$, таким же рассуждением как в пункте 1 можно добиться того, что $z$ разложится в произведение $v\cdot z'$, что, однако, уже не даёт разложимости $[z]$, т. к. $z'$ не обязано лежать в $Z$. Заметим, что условие $vu=0$ позволяет выбрать $z'$ таким образом, чтобы $u$ не содержалось в коэффициентах ни одного монома $z'$. Чтобы выделить звёздочку из $[z]$, возьмём любой моном $m$, входящий в $z'$, и пусть $uv_1$, \ldots, $uv_t$ --- все рёбра, имеющие один конец в $u$ и такие, что $m$ делится на $e_{uv_k}$. Тогда, поскольку $m$ не делится на $u$, а его степень по $u$ равна $t$, выполнено равенство $t=\deg_u z'$. Вынося в $z'$ отовсюду, откуда возможно, множитель $e_{uv_1}\wedge\cdots\wedge e_{uv_t}$, получим иное представление $z$ в виде:
$$
z=v\cdot z'=v\cdot e_{uv_1}\wedge\cdots\wedge e_{uv_t}\wedge\alpha + v\beta.
$$
Полученное разложение обладает тем свойством, что $\deg_u\alpha=0$ --- это следует из того, что $t=\deg_u z'=t+\deg_u\alpha$, и кроме того, ни один моном $\beta$ не делится на $e_{uv_1}\wedge\cdots\wedge e_{uv_t}$. Заметим, что в нашем редуцированном комплексе
$$
[v\cdot e_{uv_1}\wedge\cdots\wedge e_{uv_t}]=\starr_u(v,v_1,\ldots,v_t),
$$ а значит, чтобы применить индукцию по количеству мономов $z$ и проверить утверждение теоремы, достаточно было бы доказать, что $d\alpha=0$. Выведем это утверждение из $dz=0$:
$$
0=dz=vdz'=d(\starr_u(v,v_1,\ldots,v_t)\wedge\alpha+v\beta)=
0+(-1)^s\cdot\starr_u(v,v_1,\ldots,v_t)\wedge d\alpha+vd\beta
$$
В обеих частях равенства приравняем коэффициенты при $e_{uv_1}\wedge\cdots\wedge e_{uv_t}$ и, учитывая, что в $d\beta$ этот коэффициент равен 0, получим равенство:
$$
0=vd\alpha.
$$
Но $\Ann v=\langle u\rangle$, и теперь пришла пора вспомнить, что $\deg_ud\alpha=\deg_u\alpha=0$, так что из этого равенства уже можно вывести требуемое: $d\alpha=0$.
}

\begin{cons}
Собственные классы графа с висячей вершиной лежат в $(\starr)$.
\end{cons}
\prove{}
{ Пусть $t_i$ соответствует висячей вершине, $[z]$~--- собственный класс. Так как $[z]$~--- собственный, $\deg_{t_i}[z]>0$, так что по лемме об ограничении относительной степени $[z]\in(\S)$.
}

\begin{cons}[об алгебре гомологий дерева]
Алгебра гомологий дерева порождена звёздочками.
\end{cons}
\prove{}
{
Так как подграфы дерева являются деревьями, с помощью индукции можно установить, что несобственные классы порождаются звёздочками. Собственные классы порождаются звёздочками по предыдущему следствию.
%
%Доказываем индукцией по количеству рёбер дерева. Любой собственный подграф дерева является деревом, к которому применимо предположение индукции, поэтому идеал несобственных гомологий порождён звёздочками.
%
%Пусть класс $[z]$ лежит в максимальном идеале $\mfm$ и является собственным, $t_i$~--- какая-нибудь висячая вершина дерева. Тогда $\deg_{t_i} z\geqslant1$ (потому что иначе, $[z]$ является несобственным), следовательно, по лемме об ограничении степени, $z\in(\starr)$. Получаем требуемое включение: $\mfm\subset(\starr)$. По лемме \ref{max-ideal-gen-lemm} о порождающих максимального идеала, она влечёт нужное утверждение.
}

Также в дальнейшем понадобится знание модульной стуктуры гомологий последовательности, то есть графа с рёбрами $t_1t_2$, $t_2t_3$, \ldots, $t_{n-1}t_n$.
\begin{lemm}[о гомологиях последовательности]\label{posl}
Пусть $K$ комплекс последовательности длины $n$ (содержащей $n-1$ ребро). Тогда $\Bbbk$-векторное пространство гомологий степени $\prod_i t_i$ в $H(K)$
\begin{itemize}
\item не более, чем одномерно, если $n\bmod3\not=2$,
\item нулевое иначе.
\end{itemize}
\end{lemm}

Сперва докажем другую лемму:
\begin{lemm}[о перемещении 2-звёздочки]\label{starvar}
Пусть класс гомологий $[z]$ есть произведение звёздочки на двух рёбрах $u_0u_1$ и $u_1u_2$ на переменную $u_3$, соединённую с $u_2$ ребром, так что $[z]=u_3\cdot\starr_{u_1}(u_0,u_2)$. Тогда $[z]$ равен произведению другой звёздочки на $u_0$: $[z]= u_0\cdot\starr_{u_2}(u_1,u_3)$.
\end{lemm}
\prove{}{
Пусть $[z']=u_0\cdot\starr_{u_2}(u_1,u_3)$. Докажем $z\sim z'$ в комплексе подграфа 
$$K(u_0u_1,u_1u_2,u_2u_3,\Bbbk[u_0,u_1,u_2,u_3]),$$
из чего будет следовать гомологичность этих элементов в исходном комплексе (содержащем, возможно, больше рёбер и больше вершин). В этом подкомплексе можно произвести редукцию по двум рёбрам: $u_0u_1$ и $u_2u_3$. Теперь легко видеть, что это редукция превратит оба класса $[z]$ и $[z']$ в одно и то же: $u_0\cdot e_{u_1u_2}\cdot u_3$, что означает их равенство.
}

\prove{ леммы \ref{posl}}{
Пусть дан $z\in Z$ цикл степени $\prod_i t_i$, не гомологичный нулю. Тогда по лемме о гомологиях дерева, 
\begin{equation}\label{prodtistar}
[z]=\prod\limits_{i\in I} t_i \bigwedge\limits_j \starr_j,
\end{equation}
для некоторого множества индексов $I$ и некоторых звёздочек $\starr_j$ на парах рёбер. Среди циклов, гомологичных $z$ выберем такой, у которого сумма индексов $I$ минимальна. 

Тогда можно показать, что в разложении (\ref{prodtistar}) этого цикла $I$ является началом натурального ряда. Пусть это не так. Тогда есть индекс $i>1$ такой, что $i\in I$ и $i-1\not\in I$. Поскольку $\deg z$ тем не менее делится на $i-1$, степень одной из звёздочек в разложении (\ref{prodtistar}) должна делиться на $i-1$, а это означает, что это есть звёздочка $\starr_{t_{i-2}}(t_{i-3},t_{i-1})$. Применяя к произведению этой звёздочки на $t_i$ лемму \ref{starvar} получаем, что $z$ гомологичен элементу с меньшей суммой индексов переменных, что даёт противоречие.

Если $z$ делится хотя бы на две переменные, то он, таким образом, делится на $t_1t_2$, а значит, гомологичен нулю по лемме \ref{reduction-example}. Значит, $I=\{1\}$, либо $I=\varnothing$. Для каждого случая, звёздочки $\starr_j$ можно выбрать только одним способом, а если $n\bmod3=2$, то такого выбора и вовсе нет.
}

\subsection{Циклы}

\begin{defin}[циклического графа]
Циклическим графом или циклом длины $n$, где $n\geqslant2$ называется граф, вершинами которого служат переменные $t_1$, \ldots, $t_n$, его рёбрам отвечают мономы $t_it_{i+1}$ при $0<i<n$, а также есть ещё одно ребро, отвечающее моному $t_nt_1$.
\end{defin}

Циклический граф длины $n$ состоит из $n$ рёбер и $n$ вершин.

\begin{lemm}[о \SH-классах циклического графа]\label{nesobstv-klass-lemm}
В циклическом графе идеал несобственных классов $\mfH$ совпадает с идеалом звёздочек $(\starr)$.
\end{lemm}
\prove{}{
Во-первых, заметим, что циклический граф не является графом звёздочки, поэтому $\mfH\supset(\starr)$. Обратное включение следует из того, что любой подграф является деревом, а максимальный идеал алгебры гомологий дерева порождается звёздочками по следствию из леммы об ограничении относительной степени.
}

\begin{lemm}[о степени собственных классов циклического графа]\label{sobst-klass-deg-lemm}
Пусть $[z]\in H$~--- однородный класс гомологий степени $m$ циклического графа $\Gamma$. Если он собственный, то $m = t_1\ldots t_n$.
\end{lemm}
\prove{}{
По лемме \ref{sobstven-lemm} о степени собственных классов, $m$ содержит все переменные. Далее, если $\deg_{t_i}z>1$ для некоторого $i$, то по лемме об ограничении относительной степени, $[z]\in(\starr)$, значит, $[z]^\#=0$.
}

%\begin{cons}[о минимальной системе циклического графа]
%Пусть система порождающих алгебры гомологий циклического графа содержит все звёздочки, а также классы $[z_1]$, \ldots, $[z_t]$. Тогда $\deg[z_i]=\prod\limits_i t_i$.
%\end{cons}
%\prove{}{
%Так как подграфы циклического графа являются деревьями, идеал несобственных классов $\mfH$ порождается звёздочками по следствию из теоремы об ограничении относительной степени. Если для некоторого $i$ выполнено $\deg[z_i]\not=\prod\limits_i t_i$, то по только что доказанной лемме, $[z_i]\in\mfH$, следовательно, система порождающих не является минимальной.
%}

\subsubsection*{Удаление ребра}
\begin{cons}[из леммы \ref{addlemma} о добавлении элемента]
Пусть граф $\Gamma'$ с комплексом $K'$ получается из графа $\Gamma$ с комплексом $K$ удалением одного ребра $t_it_j$. Тогда гомологии этих комплексов связаны точной последовательностью:
\begin{equation}\label{exactkos}
H'_i \xrightarrow{u} H_i \xrightarrow{v} (0:t_it_j)_{H'_{i-1}} \to 0.
\end{equation}
\end{cons}
Бедум называть $v$ гомоморфизмом удаления ребра.

\begin{defin}[удаления ребра класса гомологий]
Если $z$ цикл в комплексе графа $\Gamma$, то класс гомологий $v([z])\in H'$ будет называться удалением ребра из $[z]$.
\end{defin}

\begin{lemm}[о поведении мультистепени при удалении ребра]
1. Пусть $[z]\in K(\x',A)$ однородный класс гомологий. Тогда $\deg u([z])=\deg z$. 
2. Пусть $[z]\in K(\x,A)$ однородный класс гомологий. Если $v([z])\not=0$, то $\deg v([z])=\frac{\deg z}{t_it_j}$.
\end{lemm}
\prove{}
{ Следует из леммы \ref{uv-deletion-lemm}.
}

\begin{lemm}[об удалении ребра звёздочки]
Пусть дана звёздочка $\starr$ на рёбрах $\alpha\beta_1$, \ldots, $\alpha\beta_k$. Тогда удаление ребра $\alpha\beta_i$ из звёздочки $\starr$ есть $$\starr_\alpha(\beta_1,\ldots,\hat\beta_i,\ldots,\beta_k),$$ либо $$-\starr_\alpha(\beta_1,\ldots,\hat \beta_i,\ldots,\beta_k),$$ где $\hat \beta_i$ обозначает отсутствие $\beta_i$.
\end{lemm}
\prove{} {
Поменяем местами $\beta_i$ и $\beta_k$, получим звёздочку $\starr_\alpha(\beta_1,\ldots,\hat \beta_i,\ldots,\beta_k, \beta_i)$, равную $\starr$ или $-\starr$. Все мономы этой звёздочки, кроме одного, содержат $e_{\alpha\beta_i}$, значит, по лемме \ref{uv-deletion-lemm}, удаление ребра $\beta_i$ из неё равно 
$$
\left[-\sum\limits_{j\not=k} (-1)^j x_j \bigwedge\limits_{s\not=j,s\not=k} e_{\alpha\beta_s}\right].
$$
}

\begin{lemm}[2-я о гомологиях последовательности]
Пусть $K$ комплекс последовательности длины $n$ (содержащей $n-1$ ребро). Тогда $\Bbbk$-векторное пространство гомологий степени $\prod_i t_i$ в $H(K)$ одномерно, если $n\bmod3\not=2$.
\end{lemm}
\prove{} {
Было доказано ранее, что указанное векторное пространство порождается одним классом гомологий $[z]$.
Этот класс есть либо произведение непересекающихся звёздочек, либо $t_1$ умноженное на такое произведение. В любом случае, повторяем операцию удаления ребра несколько раз: из каждой звёздочки удаляем её левое ребро, оставляя одно ребро. Удаление всех этих рёбер из $[z]$ превратит любую звёздочку, входящую в его разложение, в звёздочку на одном ребре, то есть переменную. Образ $[z]$ будет равен либо $[t_3t_6t_9\cdot\ldots\cdot t_{n-3}t_n]$, либо $[t_1t_4t_7\cdot\ldots\cdot t_{n-3}t_n]$. Легко видеть, что это не граница, так как в $K_1$ в данной мультистепени вообще нет элементов. Поскольку образ $[z]$ при модульном гомоморфизме оказался ненулевым, то же можно сказать и про сам $z$.
}

\subsubsection*{Гомологии цикла}
Пусть фиксирован циклический граф $\Gamma$ длины $n$. Для исследования алгебры его гомологий производится удаление ребра $t_nt_1$.

\begin{lemm}[о сведении собственных классов циклического графа к последовательности]
Пусть $H'$ обозначает аглебру последовательности $\Gamma'$, получаемой из $\Gamma$ удалением ребра $t_nt_1$. Тогда имеет место равенство:
$$
H_i/\mfH = H_i/(\S) = (0:t_nt_1)_{(H'_{i-1})}/v((\starr)_{H_i}),
$$
в котором через $v$ обозначено отображение удаления ребра формулы (\ref{exactkos}).
\end{lemm}
\prove{}{
Следует из точной последовательности (\ref{exactkos}) и из того, что 
$$\Im u\subset(\starr)=\mfH.$$
}

\begin{lemm}
В комплексе циклического графа алгебра $H/(\S)$ либо нулевая, либо одномерна как векторное пространство. Причём при $n\not\equiv2 \pmod 3$ алгебра $H/(\S)$ заведомо нулевая.
\end{lemm}
\prove{}{
По лемме \ref{sobst-klass-deg-lemm} о степени собственных классов циклического графа получаем, что $H/(\S)$ состоит из \SS-классов элементов мультистепени $\prod t_i$. По лемме о поведении мультистепени при удалении ребра получаем отсюда, что элементы $(0:t_nt_1)_{(H'_{i-1})}/v((\starr)_{H_i})$ являются образами классов гомологий со степенью $\prod\limits_{1<i<n} t_i$. По первой лемме о гомологии последовательности, векторное пространство таких классов не более, чем одномерно одномерно, тоже можно сказать и про его фактор. Причём если $n\equiv1 \pmod 3$, то по той же лемме оно нулевое. Осталось доказать тоже при $n\equiv0 \pmod 3$. В этом случае, векторное пространство классов из $H'_{i-1}$ степени $\prod\limits_{1<i<n} t_i$ порождено элементом 
$$
[z] = [t_2\bigwedge\limits_{0<i<\frac n3-1}\starr_{3i+1}(3i,3i+2)].
$$
этот элемент есть удаление ребра из разложимого класса циклического графа, равного 
$$
\starr_1(n,2)\bigwedge\limits_{0<i<\frac n3-1}\starr_{3i+1}(3i,3i+2),
$$
поэтому $[z]$ лежит в $v((\starr)_{H_i})$. В этом случае получили также $H/(\S)=0$.
}

\subsubsection*{Кружочки}
Ранее установлено, что при $n\not\equiv2 \pmod 3$ алгебра гомологий циклического графа длины $n$ порождается звёздочками. Сейчас будет показано, что при $n\equiv2 \pmod 3$ это не так, а именно, будет представлен элемент, не порождаемый звёздочками. Этот элемент назван кружочком. Так как ранее доказано, что $H/(\S)$ имеет векторную размерность не более 1, в случае $n\equiv2 \pmod 3$ кружочек и звёздочки порождают вместе алгебру гомологий графа.

Для того, чтобы получить кружочек, произведём редукцию по набору рёбер $t_{3i-2}t_{3i-1}$ для $0<i\leqslant\frac{n+1}3$. В редуцированном комплексе положим
\begin{equation}\label{eqcirc}
z = e_{t_1t_n}\bigwedge\limits_{0<i<\frac{n+1}3}\starr_{3i}(3i-1,3i+1)=
    e_{t_1t_n}\bigwedge\limits_{0<i<\frac{n+1}3}(t_{3i-1}e_{t_{3i}t_{3i+1}}-t_{3i+1}e_{t_{3i}t_{3i-1}}).
\end{equation}
Учитывая, что $t_{3i-2}t_{3i-1}=0$, замечаем, что $z$ равен сумме двух мономов:

$$
z=e_{t_1t_n}\bigwedge\limits_{0<i<\frac{n+1}3}(t_{3i-1}e_{t_{3i}t_{3i+1}})-
  e_{t_1t_n}\bigwedge\limits_{0<i<\frac{n+1}3}(t_{3i+1}e_{t_{3i}t_{3i-1}}),
$$

один из которых делится на $t_2$, а другой --- на $t_{n-1}$, так что они оба аннулируются ребром $t_1t_n$. Это и происходит при дифференцировании $z$:
\begin{align*}
dz &= {t_1t_n}\bigwedge\limits_{0<i<\frac{n+1}3}\starr_{3i}(3i-1,3i+1) = \\
   &= t_1t_n  \left(
        \bigwedge\limits_{0<i<\frac{n+1}3}(t_{3i-1}e_{t_{3i}t_{3i+1}})-
        \bigwedge\limits_{0<i<\frac{n+1}3}(t_{3i+1}e_{t_{3i}t_{3i-1}}) 
            \right) 
   = 0
\end{align*}

\begin{defin}[кружочка]
Если $n\equiv2 \pmod 3$, то класс из $K$, редуцирующийся в класс цикла (\ref{eqcirc}), называется кружочком длины $n$, и обозначается $\bigcirc_n$.
\end{defin}

Далее в этом параграфе рассматривается только случай $n\equiv2 \pmod 3$.

Степень кружочка есть произведение всех переменных, а его гомологическая степень равна $\frac{n+1}3$.

\begin{lemm}[об удалении ребра кружочка]
При удалении ребра $t_nt_1$ кружочек переходит в элемент 
$$
\bigwedge\limits_{0<i<\frac{n+1}3}\starr_{3i}(3i-1,3i+1).
$$

При удалении любого другого ребра, кружочек также переходит в произведение звёздочек.
\end{lemm}
\prove{}{
Первое утверждение следует из леммы о явном виде добавления элемента. Второе очевидно, так как пространство гомологий нужной степени одномерно и порождается произведением звёздочек.
}

\begin{lemm}[о необходимости кружочка]
$$
\bigcirc_n^\#\ne0,
$$
то есть алгебра гомологий циклического графа длины $n$ не порождается звёздочками.
\end{lemm}
\prove{}{
По лемме \ref{nesobstv-klass-lemm} о \SH-классах циклического графа, достаточно доказать \S-неразложимость кружочка. Неразложимыми элементами в меньшей мультистепени являются звёздочки, так что если бы кружочек был \S-разложим, он был бы порождён произведениями звёздочек (на паре рёбер, так как других нет), и возможно, переменных. Размерность звёздочки на паре рёбер равна 1, значит, кружочек, в случае \S-разложимости, порождался бы произведениями $\frac{n+1}3$ звёздочек. Но абсолютная степень звёздочки равна 3, так что абсолютная степень этого произведения равна $n+1$, что противоречит тому, что степень кружочка равна $n$.
}

Из сказанного выше вытекает следующий результат:
\begin{theor*}
Минимальная система порождающих для алгебры гомологий циклического графа длины $n$
\begin{itemize}
\item состоит из звёздочек на парах рёбер в случае, если $n\not\equiv2\pmod3$,
\item состоит из звёздочек на парах рёбер и кружочка длины $n$ в случае, если $n\equiv2\pmod3$.
\end{itemize}
\end{theor*}

\subsubsection*{Гомологии полной мультистепени}
Для дальнейшего понадобиться знание гомологий степени $\prod_i t_i$ циклического графа. Достаточно одного результата, который формулируется в следующей лемме:
\begin{lemm}\label{cycledimension}
Векторное пространство гомологий степени $\prod_i t_i$ циклического графа длины $n$
\begin{itemize}
\item одномерно, если $n \equiv 1$ или $n \equiv 2 \pmod 3$;
\item двумерно, если $n \equiv 0 \pmod 3$.
\end{itemize}
\end{lemm}
\prove{}{
\par\noindent \textit{Случай} $n \equiv 0 \pmod 3$. 
В этом случае ненулевым классом гомологий будет произведение $\frac n3$ непересекающихся звёздочек:
$$
[z_1]=\bigwedge\limits_{1\leqslant i\leqslant\frac n3}\starr_{3i-1}(3i-2,3i).
$$ 
Неравенство нулю доказывается стандартным способом удаления рёбер, при котором от каждой звёздочки остаётся лишь одна вершина (достаточно удалить рёбра $t_{3i-1}t_{3i}$). Кроме того, данный класс можно повернуть по часовой стрелке, получив новый класс гомологий $[z_2]$, линейно независимый с исходным. Это становится очевидно, если из их гипотетической линейной зависимости удалить ребро $t_nt_1$, при чём $[z_1]$ обнулится, а $[z_2]$ --- нет. Точно также, поворотом $[z_2]$ можно получить ещё один класс, $[z_3]$, некратный как $[z_2]$, так и $[z_1]$. Все три, однако, они линейно зависимы. В этом можно просто убедиться редуцировав по каждому третьему ребру $t_{3i}t_{3i+1}$. Классы, которые являются произведением звёздочек и переменных равны нулю, поскольку лемма \ref{starvar} позволяет приблизить эти переменные друг к другу, что будет означать делимость цикла на $t_it_{i+1}$, и следовательно, гомологичность нулю (лемма \ref{reduction-example}).

\noindent \textit{Случай} $n \equiv 1 \pmod 3$. все классы гомологий цикла данной длины порождены произведениями звёздочек на мономы. Мономы не могут содержать более одной переменной (леммы \ref{starvar}, \ref{reduction-example}), следовательно содержат одну (делимость на 3). Рассмотрим такой класс
$$
[z]=t_n\bigwedge\limits_{1\leqslant i\leqslant\frac {n-1}3}\starr_{3i-1}(3i-2,3i).
$$
Применяя к нему $\frac {n-1}3$ раз лемму \ref{starvar} получим поворот этого цикла на одну позицию, из чего следует, что все классы гомологий данной степени пропорциональны. Удаляя рёбра, домножая $[z]$ на любую переменную, непременно из него будет получаться ноль. Это может породить подозрение о том, что $z$ сам гомологичен нулю. Однако, это не так. Для того, чтобы это установить, предположим противное и представим $z$ как $dz^*$, а в $z^*$ вынесем $e_{t_nt_1}$ из всех мономов, в которых этот множитель присутствует:
$$
z=d(e_{t_nt_1}\wedge z_1^*+z_2^*).
$$
Заметим также, что $z_1^*\ne 0$, потому что иначе $z$ был бы гомологичен нулю и в комплексе последовательности рёбер, что разумеется неверно (доказывается удалением всех рёбер вида $t_{3i-2}t_{3i-1}$). А вот $z_1^*$ уже имеет степень $\prod_{2\leqslant i\leqslant n-1} t_i$, то есть является циклом последовательности длины $n-1 \equiv 2 \pmod 3$. Как известно, все такие циклы гомологичны нулю, то есть $z_1^*=dz_1^{**}$. Прибавляя к $z_1^*$ границу, не меняющую дифференциал, получаем:
$$
z=d(t_nt_1\wedge z_1^{**}+z_2^*),
$$
а значит, $z$ всё-таки получился дифференциалом элемента без $e_{t_nt_1}$, что, как уже было сказано, невозможно.

\noindent \textit{Случай} $n \equiv 2 \pmod 3$. Для того, чтобы получить разложимый класс в данной степени, понадобилось бы звёздочки умножить по крайней мере на две переменные, что даст, непременно, нулевой класс гомологий. Следовательно, все классы нужной степени порождаются именно единственным кружочком (как было доказано, все неразложимые классы равны с точностью до разложимых, которые в данной степени отсутствуют).
}
 
\subsection{Восьмёрки. Теорема о перешейке}
Начинаем параграф с полезной леммы, которая отвечает на вопрос, как выглядит аннулятор переменной $t_i$ в $H(\Gamma)$, для произвольного графа $\Gamma$, и одновременно даёт простую характеризацию идеала вершины.

\begin{lemm}[об аннуляторе переменной]\label{ann-lemma}
Пусть $t_i$ вершина графа $\Gamma$. Тогда $(0:t_i)_{H(\Gamma)}=(\starr_{t_i})$.
\end{lemm}
\prove{}{
Итак, нас интересует $(0:t_i)_H$. Добавим к нашему кольцу ещё одну переменную $t_0$, которая будет вершиной, не соединённой ни с чем ребром. Следовательно, она не может ничего аннулировать. А это означает, что $(0:t_it_0)_{H[t_0]}=(0:t_i)_{H[t_0]}=(0:t_i)_H[t_0]$. Теперь рассмотрим граф $\Gamma'$, получаемый из $\Gamma$ добавлением одного ребра $t_it_0$. Для него напишем точную последовательность удаления ребра $t_it_0$, в которую будет входить искомый аннулятор:
$$
H_i[t_0] \xrightarrow{u} H'_i[t_0] \xrightarrow{v} (0:t_it_0)_{H_{i-1}[t_0]} \to 0.
$$
Так как отображение $v$ сюрьективно, для любого цикла $z\in Z$ можно найти $\tilde z$ такой, что $[z]=v([\tilde z])$. Нас интересуют только те циклы $z$, которые не содержат $t_0$, а значит, $\deg_{t_0}(\tilde z)=1$. Лемма об ограничении степени тогда сразу же утверждает, что 
$$
[\tilde z] = \sum_s \starr_{t_i}(t_0,t_{s1},t_{s2},\ldots,t_{sk_s}) \wedge [z_s],
$$
где $z_s\in Z$ (потому, что $z_s\in Z'$ и $\deg_{t_0} z=0$).

Удалить ребро из такого элемента очень просто: нужно в каждом слагаемом у звёздочки исключить ребро $t_it_0$. Получится вот что:
$$
[z] = \sum_s \starr_{t_i}(t_{s1},t_{s2},\ldots,t_{sk_s}) \wedge [z_s].
$$

Очевидно также, что любой элемент такого вида лежит в аннуляторе.
}

\begin{defin}[восьмёрки, её звеньев]
Называем $n,m$-восьмёркой связный граф, который получается из двух циклических графов длин $n$ и $m$ склеиванием их по паре вершин. Соответствующие циклические графы называются звеньями восьмёрки.
\end{defin}

%Чтобы в достаточной мере исследовать гомологии Козюля восьмёрок, нам понадобится ещё один мощный инструмент, называемый расклейкой, принцип действия которого очень похож на удаление ребра. Но в сущности, операция расклейки несколько сложнее и требует, помимо некоторой подготовки, обращения к доказанным ранее общим фактам о комплексах над полем.

\begin{defin}[перешейка]
Пусть дан граф $\Gamma$. Его ребро называется перешейком, если при удалении этого ребра увеличивается количество компонент связности графа.
\end{defin}

Два результата этого параграфа, дающие сведенья о восьмёрках и о графах с перешейком опираются на доказанную ранее лемму о границах (лемма \ref{boardlemma}).

\begin{theor}[о перешейке]
Пусть граф $\Gamma$ при удалении ребра $\alpha\beta$ распадается на
два несвязанных графа $\Gamma_\alpha$ и $\Gamma_\beta$ с комплексами $K^\alpha=K(\Gamma_\alpha)$ и $K^\beta=K(\Gamma_\beta)$. В этом случае алгебра его гомологий $H(\Gamma)$ порождается:
\begin{itemize}
\item классами гомологий подграфов $\Gamma_\alpha$ и $\Gamma_\beta$,
\item звёздочками в вершинах $\alpha$ ($
  \starr_\alpha=\dfrac{d\bigwedge\limits_{i\in I} e_{t_i\alpha}}{\alpha}$ для некоторого множества индексов $I$) и $\beta$.
\end{itemize}
\end{theor}
\prove{}{
Прежде всего, обозначим через $\bar K$ редукцию комплекса $K=K(\Gamma)$, производимую по ребру $\alpha\beta$. Тогда $\bar H=H$ и будем доказывать утверждение теоремы для $\bar H$.

Теорема о перешейке будет доказана в два этапа:
\par\indent   1) Сперва убедимся, что $\bar H$ порождается образом $d^{-1}(\alpha Z^\alpha\otimes \beta Z^\beta)$, где $d$ --- есть дифференциал комплекса $K^\alpha\otimes K^\beta\subset K(\Gamma)$.
\par\noindent Затем в лемме \ref{factorlemma} положим $K=K^\alpha$, $C=K^\beta$, $K'=\alpha Z^\alpha$ и $C'=\beta Z^{\beta}$ и из неё получим
$$
d^{-1}(\alpha Z^\alpha\otimes \beta Z^\beta)=
d^{-1}(\alpha Z^\alpha)\otimes \beta Z^\beta + \alpha Z^\alpha\otimes d^{-1}(\beta Z^\beta) + Z(K^\alpha \otimes K^\beta).
$$
И после этого останется доказать вторую часть:
\par\indent   2) Образ подкомплекса
\begin{equation}\label{dminus1}
d^{-1}(\alpha Z^\alpha)\otimes \beta Z^\beta + \alpha Z^\alpha\otimes d^{-1}(\beta Z^\beta) + Z(K^\alpha \otimes K^\beta)
\end{equation}
в $\bar H$ порождается элементами, указанными в условии теоремы.

Итак, приступим к доказательству 1), а именно, докажем, что в каждом классе $[z] \in \bar H$ (при факторизации вначале по $\alpha\beta$, а затем по $\bar B$) можно выбрать такого представителя $z \in K^\alpha\otimes K^\beta$, что $dz\in \alpha Z^\alpha\otimes \beta Z^\beta$. Так как комплекс $\bar K$ не содержит $e_{\alpha\beta}$, то в данном классе $[z]\in\bar H$ можно выбрать представителя $z'$, лежащего в $K^\alpha\otimes K^\beta$. Тот факт, что редукция $dz'$ равна $0$ свидетельствует о том, что $dz'$ делится на $\alpha\beta$. Запишем это в виде $dz'=\alpha\beta\hat z$. Тогда $0=ddz'=\alpha\beta d\hat z$, а значит, $d\hat z=0$, потому что это равенство выполнено уже в нередуцированном комплексе. Применяя к $\hat z$ формулу Кюннета в виде (\ref{kunetform}), получим его представление в виде $\hat z=\hat z_1 + \hat z_2$ с $\hat z_1\in Z^\alpha\otimes Z^\beta$ и $\hat z_2\in B$. Теперь выберем элемент $\hat z_2'$, такой что $d\hat z_2'=\hat z_2$ и положим $z=z'-\alpha\beta\hat z_2'$. Редукции $[z]$ и $[z']$ совпадают, т. к. они отличаются на элемент кратный $\alpha\beta$, но дифференциал $dz$ можно записать как
$$
dz=d(z'-\alpha\beta\hat z_2')=\alpha\beta\cdot(\hat z_1 + \hat z_2)-\alpha\beta\hat z_2=\alpha\beta\hat z_1,
$$
про $\hat z_1$ же известно, что он лежит в $Z^\alpha\otimes Z^\beta$, так что 1) доказано.

Доказательство 2). Докажем этот факт отдельно для каждого из трёх слагаемых (\ref{dminus1}). Для $Z(K^\alpha\otimes K^\beta)$ он прямо следует из формулы Кюннета. Для $d^{-1}(\alpha Z^\alpha)\otimes \beta Z^\beta$ (и совершенно аналогично для $\alpha Z^\alpha\otimes d^{-1}(\beta Z^\beta)$) получается несколько сложнее. А именно, выберем однородный $z\in d^{-1}(\alpha Z^\alpha)\otimes \beta Z^\beta$, и будем доказывать  индукцией по количеству мономов в записи $z$, не делящихся на $\alpha$. Если таковых мономов нет, то $z$ само делится на $\alpha$ и значит на $\alpha\beta$, что означает гомологичность нулю. Рассмотрим максимальное по длине произведение $\bigwedge\limits_i e_{x_i\alpha}$, участвующее в записи $z$. Однородность тогда означает, что там будет столько множителей $e_{x_i\alpha}$, какова относительная степень $\deg_\alpha z$. Вынесем произведение $\bigwedge\limits_i e_{x_i\alpha}$ за скобки во всех мономах $z$, в которых оно есть и получим новое представление 
$$
z=\bigwedge\limits_i e_{x_i\alpha}\wedge\beta z_1 + z_2,
$$
где $\deg_\alpha z_1=0$, а $z_2$ не содержит произведения $\bigwedge\limits_i e_{x_i\alpha}$. Обобщённая редукция $\bar z$ элемента $z$ есть цикл в $\bar K$, так что $\beta\bar z_1$ есть цикл (это получается как обычно приравниванием к нулю коэффициента при $\bigwedge\limits_i e_{x_i\alpha}$ в $d\bar z$). Про $\beta z_1$ это обстоятельство означает, что его дифференциал делится на $\alpha\beta$, но $\deg_\alpha\beta dz_1=\deg_\alpha\beta z_1=0$, откуда уже следует, что эта делимость возможна только в случае $dz_1=0$. Формула Кюннета опять порождает этот элемент нужными по условию, а класс обобщённой редукции $e_{x_i\alpha}\cdot\beta$ есть не что иное, как звёздочка $\starr_\alpha(\beta,x_1,x_2,\ldots)$. Для второго слагаемого $z_2$ в такой ситуации ничего не остаётся, как самому лежать в $d^{-1}(\alpha Z^\alpha\otimes \beta Z^\beta) + \bar B$, а значит, к нему применима индукция, которая и доказывает теорему окончательно.
}

\subsubsection*{Склейка и расклейка}
\begin{defin}[склейки и расклейки графов]
Пусть допустимый (без петель и кратных рёбер) граф $\Gamma'$ с комплексом $K(\Gamma')$ получается из допустимого графа $\Gamma$ с комплексом $K(\Gamma)$ отождествлением его вершин $\alpha$ и $\beta$ в одну вершину $\gamma$. Тогда граф $\Gamma'$ называется склейкой $\Gamma$ в вершинах $\alpha$ и $\beta$. Граф $\Gamma$ называется расклейкой $\Gamma'$ в вершине $\gamma$.
\end{defin}

Для того, чтобы при отождествлении вершин получался граф, удовлетворяющий требуемым условиям для построения комплекса Козюля, нужно наложить два ограничения на его рёбра:
\begin{itemize}
\item склеенные вершины не соединены ребром (отсутствие петель);
\item не существует вершины, соединённой рёбрами как с $\alpha$, так и с $\beta$ (отсутствие кратных рёбер).
\end{itemize}

Расклейка определяется не только расклеиваемой вершиной, но также тем, какие рёбра, их соединённых в склеенном графе с вершиной $\gamma$ будут после расклейки соединены с $\alpha$, а какие с $\beta$.

\begin{lemm}[о точной последовательности расклейки]
Пусть $\Gamma'$ является склейкой $\Gamma$ в вершинах $\alpha$ и $\beta$.

Имеет место точная последовательность модулей:
\begin{equation}\label{raskleika-exact}
H_i(\Gamma) \xrightarrow{u} H_i(\Gamma') \xrightarrow{v} (0:\alpha-\beta)_{H_{i-1}(\Gamma)} \to 0
\end{equation}
\end{lemm}
\prove{}{Посмотрим теперь на комплекс $K(\Gamma')$. Он получается из комплекса $K(\Gamma)$ факторизацией по $(\alpha-\beta)$, то есть на самом деле равен $K(\Gamma,A/(\alpha-\beta))$. После того, как мы записали комплекс в таком виде, в нему можно применить лемму о редукции, получая такой изоморфизм:
$$
H(\Gamma')\cong H(\Gamma,\alpha-\beta,A).
$$
В правой части, однако, стоит уже не мономиальный комплекс, так что большинство доказанных фактов не будет для него работать. Однако, лемму  о добавлении элемента применить можно, и получить точную последовательность (\ref{exact}), которая в данном случае принимает вид:
$$
H_i(\Gamma) \xrightarrow{u} H_i(\Gamma') \xrightarrow{v} (0:\alpha-\beta)_{H_{i-1}(\Gamma)} \to 0.
$$
}

\begin{defin}[расклейки класса гомологий]
Образ класса $[z]\in H_i(\Gamma')$ при отображении $v$ называется \textit{расклейкой класса гомологий $[z]$}.
\end{defin}

\begin{lemm}[явный вид отображений $u$ и $v$ при расклейке]
Отображение $u$ является отображением факторизации. $v$ переводит класс $[z]$ в элемент $\frac{dv_1([z])}{\alpha-\beta}$, где 
$$
H_i(\Gamma') \xrightarrow{v_1} 
H_i(\Gamma,A/(\alpha-\beta)).
$$
\end{lemm}
\prove{}{
Утверждение об $u$ следует из его определения.

Отображение $v$ представим в виде композиции следующим образом:
\begin{equation}\label{vexpand}
v:
H_i(\Gamma') \xrightarrow{v_1} 
H_i(\Gamma,A/(\alpha-\beta)) \xrightarrow{v_2} 
H_i(\Gamma,\alpha-\beta,A) \xrightarrow{v_3}
H_{i-1}(\Gamma).
\end{equation}

Теперь опишем, как преобразуется класс гомологий $[z]$, представленный данным циклом $z$, при этих отображениях:
\begin{itemize}
\item $v_1([z])$ есть класс цикла, получающегося из $z$ заменой каждого вхождения $e_{\gamma t_i}$ на $e_{\alpha t_i}$ или на $e_{\beta t_i}$, в зависимости от того, какое ребро ($\alpha t_i$ или $\beta t_i$) есть в расклеенном графе, а также заменой каждого вхождения $\gamma$ произвольно: на $\alpha$ или $\beta$ (поскольку образы $\alpha$ и $\beta$ в кольце $A/(\alpha)$ они равны).
\item $v_2$ есть поднятие при редукции, то есть, согласно лемме \ref{podniatie}, 
\begin{equation}\label{v2v1}
v_2(v_1([z]))=v_1([z])-e_{\alpha-\beta}\wedge\frac{dv_1([z])}{\alpha-\beta}.
\end{equation}
\item Наконец, $v_3$ есть отображение добавления элемента, участвующее в (\ref{exact}). Оно (\ref{v2v1}) переводит в 
$\frac{dv_1([z])}{\alpha-\beta}$.
\end{itemize}
}

\begin{cons}[о расклейке произведения]\label{prodimagelemma}\label{raskleim-prod}
Для отображения $v$ справедливо равенство (\ref{prodimage}):
$$
u\circ v([z]\wedge [z'])=u\circ v([z])\wedge[z']-(-1)^{\dim z}[z]\wedge u\circ v([z']).
$$
\end{cons}
\prove{}{
Поскольку $v_1$ и $v_2$ есть отображения алгебр, с произведением $\wedge$ они перестановочны, а к $v_3$ можно применить лемму \ref{prodimagelemma}, и получить требуемое.
}

В общем случае, рсклейка однородного класса гомологий $z$ не является однородным элементом. Нарушение однородности в разложении (\ref{vexpand}) происходит только при $v_1$ и затрагивает только отождествляемые переменные $\alpha$ и $\beta$.

\begin{lemm}[о нарушении однородности]
Пусть $[z]$ однородный класс гомологий, а $v$~--- отображение, расклеивающее вершину $\gamma$ в $\alpha$ и $\beta$. Тогда для лобого монома $m$ из $v(z)$ выполнено соотношение
$$
\deg_\alpha m+\deg_\beta m=\deg_\gamma [z]-1,
$$
$$
\deg_t m=\deg_t [z], \mbox{ если вершина $t$ отлична от $\alpha$ и $\beta$.}
$$

Таким образом, однородные компоненты $v([z])$ имеют степени 
$$\frac{\deg z}{\gamma^{\deg_\gamma z}}\alpha^{r_\alpha}\beta^{r_\beta}, \mbox{ где $r_\alpha+r_\beta=\deg_{\gamma} z-1$.}$$
\end{lemm}
\prove{}{Очевидно следует из явного вида $v$.}

Обозначим через $v_{r_\alpha r_\beta}([z])$, где $r_\alpha$ и $r_\beta$~--- неотрицательные целые числа, однородную компоненту $v([z])$, имеющую степень
$$
\frac{\deg z}{\gamma^{\deg_\gamma z}}\alpha^{r_\alpha}\beta^{r_\beta}.
$$

\begin{lemm}[ограничение вида однородных компонент расклейки]\label{raskleika-comp-lemma}
Пусть валентность вершины $\gamma$ равна $s$, в расклеенном графе из вершины $\alpha$ выходит $s_\alpha$ рёбер, а из вершины $\beta$ --- $s_\beta$ рёбер (тогда $s=s_\alpha+s_\beta$). Тогда если $v_{r_\alpha r_\beta}([z])\not=0$, то $r_\alpha<s_\alpha$ и $r_\beta<s_\beta$.
\end{lemm}
\prove{}{ Предположим, например, что однородная компонента $v_{r_\alpha r_\beta}([z])$ имеет максимальный индекс $r_\alpha$ среди однородных компонент $v([z])$, причём этот индекс не меньше $s_\alpha$. Вспомним, что $v([z])\in(0:\alpha-\beta)$, так что $v([z])\cdot(\alpha-\beta)=0$. Выделим в произведении $v(z)\cdot(\alpha-\beta)=0$ однородную компоненту, содержащую $\alpha ^ {r_\alpha+1}$. Эта компонента гомологична $\alpha v_{r_\alpha r_\beta}(z)$, так как у всех других слагаемых относительная степень по $\alpha$ меньше указанного числа. Таким образом, в $K$ выполнено $\alpha v_{r_\alpha r_\beta}(z)=dz^*$ для некоторого $z^* \in K$. Поскольку относительная степень $z^*$ равна $r_\alpha+1$, что больше степени вершины $\alpha$, $z^*$ делится на $\alpha$, откуда следует, что $\alpha v_{r_\alpha r_\beta}(z)=\alpha d(\frac{z^*}{\alpha})$. Обе части можно сократить на $\alpha$, потому что равенство выполнено в комплексе $K$ без делителей нуля, так что в $H$ имеет место равенство: $v_{r_\alpha r_\beta}([z])=0$.
}
\begin{cons}[о максимальном количестве однородных компонент, возникающих при расклейке]\label{maxcompcconst}
В $v([z])$ не может быть однородных компонент больше, чем $\min(\deg\alpha,\deg\beta)$. 
\end{cons}
\begin{lemm}[о связи между однородными компонентами расклейки]\label{raskleika-components-link}
Для произвольного натурального $k$ имеем:
\begin{equation}
\beta^k v_{r_\alpha+k,r_\beta}([z]) = \alpha^k v_{r_\alpha,r_\beta+k}([z]).
\end{equation}
\end{lemm}
\prove{}{
Рассмотрим тождество 
$$
v([z])(\alpha-\beta)=0,
$$
из которого следует равенство нулю его однородных компонент. Запишем это:
$$
\beta v_{r_\alpha+1,r_\beta}([z]) = \alpha v_{r_\alpha,r_\beta+1}([z]).
$$
Применяя это равенство несколько раз можно получить утверждение леммы.
}

\begin{lemm}[о расклейке звёздочки]\label{raskleim-starr}
Пусть $\Gamma'$ поулчается из $\Gamma$, склеиванием вершин $\alpha$ и $\beta$ в вершину $\gamma$, в графе $\Gamma$ имеются рёбра $\alpha r_1$, \ldots, $\alpha r_p$ и $\beta s_1$, \ldots, $\beta s_q$. Тогда расклейка звёздочки комплекса $\Gamma'$
$$
\starr_\gamma(\r,\s)
$$
равна произведению двух звёздочек со знаком $+$ или $-$:
$$
\pm\starr_\alpha(\r)\wedge\starr_\beta(\s)
$$
\end{lemm}
\prove{}{Прямое вычисление, которое значительно упрощается, если редуцировать по ребру звёздочки (редукция по ребру корректна, потому что последовательность $\alpha-\beta,\alpha t_i$ регулярна), показывает, что утверждение верно. Проверим это явно.%расклейка звёздочки равна просто-напросто произведению двух звёздочек, одна из которых находится в вершине $\alpha$, другая --- в вершине $\beta$.

Произведём редукцию по ребру $\alpha r_1$ (в редуцированном комплексе --- по ребру $\gamma r_1$), тогда звёздочка $\starr_\gamma(\r,\s)$ принимает вид
$$
\bar\starr=\left[r_1\bigwedge_{i\not=1}e_{r_i\gamma}\bigwedge_j e_{s_j\gamma}\right].
$$

Чтобы вычислить $v(\bar\starr)$ воспользуемся разложением (\ref{vexpand}). Имеем:
$$
v_1(\bar\starr)=\left[r_1\bigwedge_{i\not=1}e_{r_i\alpha}\bigwedge_j e_{s_j\beta}\right].
$$

Чтобы вычислить образ $v_2$, требуется вычисление $dv_1(\bar\starr)$, которое упрощается, если заметить, что $d(r_1\bigwedge_{i\not=1}e_{r_i\alpha})=0$:
$$
dv_1(\bar\starr)=\pm\left[r_1\bigwedge_{i\not=1}e_{r_i\alpha} \wedge d\bigwedge_j e_{s_j\beta}\right]=
\pm\left[r_1\bigwedge_{i\not=1}e_{r_i\alpha}\right]\wedge \beta\starr_{\beta}(\s).
$$

Поскольку $r_1\alpha=0$ выполнено $r_1\beta=-r_1(\alpha-\beta)$, так что полученное представление $dv_1(\bar\starr)$ можно ещё преобразовать:
$$
dv_1(\bar\starr)=
\pm(\alpha-\beta)\starr_{\alpha}(\r)\wedge\starr_{\beta}(\s),
$$
откуда уже следует, что $v_2(v_1(\bar\starr))$ равно
$$
v_1(\bar\starr)\pm e_{\alpha-\beta}\wedge\starr_{\alpha}(\r)\wedge\starr_{\beta}(\s).
$$
Наконец, осталось применить отображение $v_3$, которое отбрасывает все слагаемые без $e_{\alpha-\beta}$, а у остальных отбрасывает $e_{\alpha-\beta}$. Получится утверждение леммы.
}

\begin{defin}[распада при расклейке]
Будем говорить, что при расклейке вершины $\gamma$ в $\alpha$ и $\beta$ имеет место распад в том случае, если в расклеянном графе вершины $\alpha$ и $\beta$ находятся в разных компонентах связности. 
\end{defin}

Когда имеет место распад при расклейке вершины $\gamma$ в $\alpha$ и $\beta$, будем разбивать граф на две несвязные  компоненты, одна из которых содержит $\alpha$, а другая~--- $\beta$. Эти компонент обозначим, соотетственно, $\Gamma_\alpha$ и $\Gamma_\beta$, их комплексы~--- $K_\alpha$ и $K_\beta$, а их гомологии~--- $H_\alpha$ и $H_\beta$.

\begin{lemm}[об аннуляторе $(\alpha,\beta)$ в случае распада]\label{annalphabeta-lemm} Пусть вершина $\gamma$ графа $\Gamma'$ расклеивается в вершины $\alpha$ и $\beta$ графа $\Gamma$, причём имеет место распад.

Тогда подмодуль $H(\Gamma)$ классов гомологий, аннулируемых одновременно и $\alpha$ и $\beta$, порождается расклейками классов вида
$$
\starr_\gamma(\k)\wedge [\tilde z],
$$
для некоторого цикла $\tilde z\in Z(\Gamma)$ и некоторых вершин $k_1$, \ldots, $k_p$, соединённых c $\gamma$ в $\Gamma'$ и обозначенных как $\k$, то есть содержится в $v((\starr_\gamma))$.
\end{lemm}

\prove{}{
Предположим, расклеив некоторый класс гомологий $[z]\in H'$, мы получили класс $v([z])\in H$, который обладает тем свойством, что $\alpha v([z])=\beta v([z])=0$. Пользуясь изоморфизмом 
\begin{equation}\label{capeq}
M'\otimes N\cap M \otimes N'=M'\otimes N',
\end{equation}
существующим для любых векторных пространств $M$ и $N$, и их подпространств, соответственно, $M'$ и $N'$, получаем 
\begin{equation}\label{annalphabeta}
v([z])\in (0:\alpha) \otimes H_\beta \cap H_\alpha \otimes (0:\beta) = (0:\alpha)_{(H_\alpha)} \otimes (0:\beta)_{(H_\beta)}.
\end{equation}

Тогда учитывая результат леммы об аннуляторе переменной, получаем, что $v(z)$ представляется в виде суммы слагаемых вида 
$$
\starr_\alpha (t_1,\ldots,t_k) \wedge\starr_\beta (s_1,\ldots,s_l) \wedge [z_1]\wedge [z_2],
$$
в которых $z_1\in Z_\alpha$, $z_2\in Z_\beta$, $t_i$ вершины $\Gamma_\alpha$, а $s_i$ вершины $\Gamma_\beta$. Но тогда, такое слагаемое является расклейкой класса гомологий $\Gamma'$, равного 
$$
\starr_\gamma (t_1,\ldots,t_k,s_1,\ldots,s_l) \wedge [z_1]\wedge [z_2].
$$
}

\begin{cons}[об однородной расклейке в случае распада]\label{raskleikaodnorodna}
Пусть имеет место распад при расклейке $\Gamma'$ в $\Gamma$, $[z]\in H(\Gamma')$. Тогда если $v([z])$ однородный элемент, то он является суммой элементов вида (\ref{annalphabeta}).
\end{cons}
\prove{}{
Заметим, что если $v([z])$ однородный элемент, то из уравнения $v([z])(\alpha-\beta)=0$, которому он удовлетворяет, можно вывести $\alpha v([z])=\beta v([z])=0$, что, как только что выяснено, даёт требуемое представление.
}

\begin{cons}[из лемм \ref{raskleika-comp-lemma}, \ref{raskleika-components-link} о соотношениях между однородными компонентами в случае $\deg\alpha=2$]
Пусть вершина $\gamma$ графа $\Gamma'$ расклеивается в вершины $\alpha$ и $\beta$ графа $\Gamma$ и имеет место распад. Если $\deg\alpha=2$, и $[z]\in H(\Gamma')$~---однородный класс гомологий, то $v([z])$ имеет не более двух ненулевых однородных компонент: $v_\alpha([z])$ и $v_\beta([z])$, которые обозначим $[z_\alpha]$ и $[z_\beta]$. Эти компоненты удовлетворяют соотношениям:
\begin{align}
\alpha [z_\alpha] &= 0, \label{firstc} \\
\beta  [z_\beta]  &= 0, \label{second} \\
\alpha [z_\beta]  &= \beta [z_\alpha]. \label{third}
\end{align}
\end{cons}

\begin{lemm}[о \SS-классах при распаде в случае $\deg\alpha=2$]\label{sclass-lemma}
Пусть вершина $\gamma$ графа $\Gamma'$ расклеивается в вершины $\alpha$ и $\beta$ графа $\Gamma$, имеет место распад и $\deg\alpha=2$. %Тогда образ множества 
%$$
%\{ [z] \in H(\Gamma') \;|\; v([z]) = \alpha[\tilde z] + \beta[\tilde z]\; \mbox{для некоторого $[\tilde z] \in H(\Gamma)$} \}
%$$
%при факторизации по $(\S)$ совпадает со всей алгеброй $H(\Gamma')/(\S)$.
Тогда в любом \SS-классе $H(\Gamma')$ есть элемент, расклейка которого имеет вид 
\begin{equation}\label{sclass-form}
\alpha[\tilde z] + \beta[\tilde z].
\end{equation}
\end{lemm}
\prove{}{
%Докажем, что в любом \SS-классе $H(\Gamma')$ есть элемент, расклейка которого имеет вид $\alpha[\tilde z] + \beta[\tilde z]$. 

Посмотрим на $(\ref{third})$. Обе части равенства лежат в 
$$
\alpha H_\alpha \otimes H_\beta \cap H_\alpha \otimes \beta H_\beta = \alpha H_\alpha \otimes \beta H_\beta
$$
(равенство векторных пространств в соответствии с изоморфизмом (\ref{capeq})). Значит, для некоторого класса гомологий $[\hat z]\in H$ выполнено
\begin{equation}\label{hatzcond}
\alpha [z_\beta] = \beta [z_\alpha] = \alpha\beta[\hat z].
\end{equation}

Рассмотрим сперва одну часть этого равенства, ту, которая утверждает, что $\alpha [z_\beta]=\alpha\beta[\hat z]$, то есть $\alpha([z_\beta] - \beta[\hat z])=0$. Это равенство домножим на $\beta$ и, воспользовавшись (\ref{second}), получим $\alpha\beta^2[\hat z]=0$. Теперь поднимем это равенство из алгебры гомологий в комплекс, получая включение $\alpha\beta^2\hat z \in \alpha Z_\alpha \otimes \beta^2 Z_\beta \cap B(K_\alpha\otimes K_\beta)$. Это пересечение можно преобразовать в соответствии с леммой \ref{boardlemma}, представляя $\hat z$ в виде суммы двух циклов: $\hat z=\hat z_\alpha+\hat z_{\beta^2}$, таких что $\hat z_\alpha\in (0:\alpha)_H$, $\hat z_{\beta^2}\in (0:\beta^2)_H$. Без нарушения условия (\ref{hatzcond}) можно заменить $\hat z$ на $\hat z_{\beta^2}$, добившись для нового $\hat z$ равенства $\beta^2[\hat z]=0$. Заметим, что если уже было выполнено $\alpha^2[\hat z]=0$ для старого $\hat z$, то это останется и для нового, потому что они отличаются на элемент, аннулируемый $\alpha$. Это означает, что повторяя те же рассуждения, поменяв только местами $\alpha$ и $\beta$, можно заменить $\hat z$ на такой элемент, который удовлетворяет (\ref{hatzcond}), а также ещё двум условиям:
\begin{align}
\alpha^2\hat z &= 0, \\
\beta ^2\hat z &= 0. \label{hatsecond}
\end{align}

Посмотрим теперь ещё раз на соотношение $\alpha[z_\beta - \beta\hat z]=0$. Оно означает, что $[z_\beta - \beta\hat z] \in (0:\alpha)$. С другой стороны, оба слагаемых этого элемента аннулируются с помощью $\beta$ (следует из (\ref{second}) и (\ref{hatsecond})), тогда по лемме \ref{annalphabeta-lemm} он есть образ $\starr$-разложимого. Аналогично заключаем и про $[z_\alpha - \alpha\hat z]$. Вычитая из $z$ эти два \S-разложимых элемента, образами классов которых являются $[z_\beta - \beta\hat z]$ и $[z_\alpha - \alpha\hat z]$, получим элемент $\tilde z$, такой что $v([\tilde z])=[\alpha\hat z + \beta\hat z]$.
}

\subsubsection*{Гомологии восьмёрки}
\begin{lemm}[о несобственных классах восьмёрки]
Если дана $(n,m)$-восьмёрка $\Gamma'$, то идеал несобственных классов $\mfH$ порождается:
\begin{itemize}
\item звёздочками,
\item кружочком в первом звене, если $n\pmod3=2$,
\item кружочком во втором звене, если $m\pmod3=2$.
\end{itemize}
\end{lemm}
\prove{}{
Доказательство тривиально в силу того, что несобственный подграф звёздочки является либо деревом, либо циклическим графом с дополнительным деревом в одной из вершин.
}

Далее будем производить расклейку восьмёрки $\Gamma'$ в её центральной вершине $\gamma$, при которой расклейка распадается на два циклических графа: $\Gamma_\alpha$ и $\Gamma_\beta$. Считаем, что вершина $\gamma$ расклеивается в вершины $\alpha$ и $\beta$, лежащие, соответственно, в $\Gamma_\alpha$ и $\Gamma_\beta$.

Пусть $\Gamma'$, $\Gamma_\alpha$ и $\Gamma_\beta$ --- те же; $\alpha$, $\alpha_1$, \ldots, $\alpha_{n-1}$ --- множество вершин графа $\Gamma_\alpha$; $\beta$, $\beta_1$, \ldots, $\beta_{m-1}$ --- множество вершин графа $\Gamma_\beta$. И пусть $U_\alpha\subset H_\alpha$ обозначает векторное пространство классов гомологий графа $\Gamma_\alpha$ мультистепени, равной $\prod\limits_{1\leqslant i\leqslant n-1} \alpha_i$; а $U_\beta\subset H_\beta$ обозначает векторное пространство классов гомологий графа $\Gamma_\beta$ мультистепени, равной $\prod\limits_{1\leqslant i\leqslant m-1} \beta_i$. Пусть $U = U_\alpha\otimes U_\beta$.

Рассмотрим \SH-класс $[z]_\#$. По лемме \ref{sclass-lemma} в нём есть представитель с расклейкой вида $\alpha[\tilde z] + \beta[\tilde z]$ для некоторого $[\tilde z] \in H(\Gamma)$. Рассмотрим базис $[z_1]_\#$, \ldots, $[z_s]_\#$ векторного пространства $H'_\#$. Для каждого $[z_i]_\#$ выберем представителя с расклейкой $\alpha[\tilde z_i] + \beta[\tilde z_i]$. Теперь рассмотрим линейное отображение $\omega$, переводящее каждый \SH-класс $[z_i]_\#$ в $[\tilde z_i]$.

\begin{lemm}[о свойствах отображения $\omega$]\label{inject-lemma}
Пусть 
$$
V=U \cap (0:\alpha^2)_H \cap (0:\beta^2)_H.
$$
Тогда построенное отображение $\omega$ является инъективным отображением
$$
\omega: H'_\#\to V.
$$
\end{lemm}
\prove{}{
Сперва докажем, что $\Im\omega\subset V$. Докажем включение $[\tilde z_i]\in U$. Учитывая результат следствия \ref{raskleikaodnorodna} получаем, что $[\alpha \tilde z_i]\ne 0$ и $[\beta \tilde z_i]\ne 0$. Тогда можно воспользоваться ограничением вида однородных компонент расклейки (леммой \ref{raskleika-comp-lemma}) и получить, что $\deg_\alpha \alpha \tilde z_i < 2$, $\deg_\beta \beta \tilde z_i < 2$, значит $\deg_\alpha \tilde z_i = 0$ и $\deg_\beta \tilde z_i = 0$. Поскольку $[z'_i]$ был собственным классом, его степень относительно всех вершин, отличных от $\gamma$ должна быть равна 1 (по лемме об ограничении относительной степени). Такая же степень, следовательно, будет и у $\tilde z_i$ относительно вершин $\alpha_j$ и $\beta_j$. Таким образом, установлено, $[\tilde z_i]\in U$.

Докажем, что $[\tilde z_i]\in (0:\alpha^2)_H$. Для этого достаточно воспользоваться (\ref{firstc}), с учётом того, что в нашем случае $[z_\alpha]=\alpha[\tilde z_i]$. Аналогично проверяется $[\tilde z_i]\in (0:\beta^2)_H$. Таким образом, установлено включение $[\tilde z_i]\in V$, следовательно, отображение
$$
\omega:[z'_i]_\# \mapsto [\tilde z_i]
$$
определено корректно. Также из линейности следует, что не только для базисных элементов, а для произвольного собственного \#-класса $[z']_\#$ верно, что расклейка некоторого его представителя равна $(\alpha+\beta)\omega([z']_\#)$.

Докажем его инъективность. Пусть $\omega([z']_\#)=0$ для некоторого собственного \#-класса $[z']_\#$, причём можно считать, что $[z']$ --- это тот представитель $[z']_\#$, расклейка которого равна $(\alpha+\beta)\omega([z']_\#)=0$. Получаем, что $[z']$ лежит в ядре расклейки, а значит, также и в образе отображения $u$ последовательности (\ref{raskleika-exact}), то есть является собственным, следовательно, $[z']_\#=0$.
}

\begin{lemm}[о сохранении инъективности $\omega$ при проекции]\label{inject-factor-lemma}
Обозначим
$$
V'=V/(U \cap ((0:\alpha)_H + (0:\beta)_H)),
$$
и пусть $\varphi$~--- естественная проекция $V\to V'$.

Тогда отображение $\varphi\circ\omega: H'_\#\to V'$ является инъективным.
\end{lemm}
\prove{}
{ Предположим противное, что $\Im\omega\cap\Ker\varphi\not=0$. Пусть сперва $\omega([z']_\#)\in (0:\alpha)_H$ причём без ограничения общности можно полагать, что $[z']$ есть тот представитель $[z']_\#$, расклейка которого равна $(\alpha+\beta)\omega([z']_\#)$. По предположению, расклейка $[z']$ тогда равна $\beta\omega([z']_\#)$. Теперь применяем следствие \ref{raskleikaodnorodna}, получая \S-разложимость $[z']$ из которой следует равенство $[z']_\#=0$. Таким образом, инъективность сохраняется при проектировании на $V/(U \cap (0:\alpha)_H)$. Аналогичное рассуждение проходит и для $\beta$.
}

\begin{theor}\label{8theor}
Алгебра \#-классов $(n,m)$-восьмёрки
\begin{itemize}
\item нулевая, если $n\not\equiv1 \pmod 3$ или $m\not\equiv1 \pmod 3$;
\item одномерна, если $m\equiv n\equiv 1 \pmod 3$.
\end{itemize}
\end{theor}
\prove{} {
Пусть $n\equiv 0\pmod 3$. В этом случае пространство $U_\alpha$ изоморфно пространству гомологий полной степени последовательности из $n-1$ вершины. Это пространство нулевое согласно лемме \ref{posl} о гомологиях последовательности. Поэтому в этом случае имеет место равенство $U=0$. В силу отображения $\omega$, которое вкладывает $H_\#$ в $U$ по лемме \ref{inject-lemma}, получаем, что $H'_\#=0$ в этом случае. Аналогично рассматривается случай $m\equiv 0\pmod 3$.

Пусть $n\equiv 2\pmod 3$. По лемме \ref{posl} о гомологиях последовательности, пространство $U_\alpha$ одномерно. В этом случае, однако, единственный базисный элемент $U_\alpha$ аннулируется $\alpha$ (в соответствии с леммой об аннуляторе переменной), так что определённое выше пространство $V'$ в этом случае нулевое. Получаем, что $H'_\#=0$ в силу леммы \ref{inject-factor-lemma}, которая утверждает наличие инъективного отображения $H'_\#\to V'$. Аналогично рассматривается случай $m\equiv 0\pmod 3$.

Остался неразобранным единственный случай $m\equiv n\equiv 1 \pmod 3$. В нём пространства $U_\alpha$ и $U_\beta$ одномерны, поэтому алгебра $H'_\#$ не более, чем одномерна. Заметим, что в данном случае, модуль $(0:\alpha-\beta)_{H(\Gamma)}$ содержит элемент 
\begin{equation}\label{netrivialniyv8}
[z]=(\alpha+\beta)\bigwedge\limits_{1\leqslant i\leqslant\frac{m-1}3}\starr_{\alpha_{3i-1}}(\alpha_{3i-2},\alpha_{3i}) 
            \wedge\bigwedge\limits_{1\leqslant i\leqslant\frac{n-1}3}\starr_{\beta _{3i-1}}(\beta _{3i-2},\beta _{3i}),
\end{equation}
где $\alpha_i$, $\alpha$ --- все вершины графа $\Gamma_\alpha$, а $\beta_i$ и $\beta$ --- все вершины графа $\Gamma_\beta$. По леммам \ref{cycledimension} и \ref{ann-lemma} видно, что этот класс гомологий не аннулируется $\alpha$, в то время как расклейка любого \S-разложимого класса гомологий обязана аннулироваться как $\alpha$, так и $\beta$ (лемма \ref{raskleim-starr}). При данных длинах звеньев кружочков в графе не существует, значит, идеал собственных классов $\mfH$ порождается звёздочками и данное рассуждение является доказательством требуемого факта.
}

Рассмотрим теперь другую расклейку в центральной вершине~--- ту, при которой восьмёрка $\Gamma'$ превращается в один циклический граф $\hat\Gamma$. Это будет циклический граф длины $m+n\equiv 2 \pmod 3$, в котором существует кружочек $\hat\bigcirc$. Отображение склейки и расклейки в этом случае обозначим как $\hat u$ и $\hat v$.

\begin{lemm}
В случае $m\bmod 3=n\bmod 3=1$ минимальная система порождающих восьмёрки состоит из минимальной системы звёздочек и склейки класса $\hat\bigcirc$, то есть $\hat u(\hat\bigcirc)$.
\end{lemm}
\prove{}
{ Пусть центральная вершина восьмёрки $\gamma$ расклеилась в две вершины: $\hat\alpha$ и $\hat\beta$. При этом нижнее звено, помимо $\gamma$, имело вершины $\alpha_1$, \ldots, $\alpha_{n-1}$, а верхнее --- $\beta_1$, \ldots, $\beta_{m-1}$. Вершины пронумерованы так, чтобы во время движения по циклическому графу $\hat\Gamma$, индексы либо всё время убывали, либо всё время возрастали (в зависимости от направления), за исключением моментов перехода на $\hat\alpha$ и $\hat\beta$. Также будут обозначены соответствующие вершины графа $\hat\Gamma$. Рассмотрим редукцию комплекса $\hat K$, производимую по следующим рёбрам:
в нижнем звене редукция производится по рёбрам $\alpha_{3i-2}\alpha_{3i-1}$ для $i=1$, \ldots, $\frac{n-1}3$; в верхнем звене редукция производится по ребру $\hat\beta\beta_1$, а также по рёбрам $\beta_{3i-1}\beta_{3i}$ для $i=1$, \ldots, $\frac{m-1}3$. Образ при такой редукции кружочка $\hat\bigcirc$ принимает известный вид
\begin{equation}\label{raskleikacirc}
e_{\alpha_{n-1}\hat\beta}\starr_{\hat\alpha}(\beta_{m-1}\alpha_1)\bigwedge_{1\leqslant i<\frac{n-1}3}\starr_{\alpha_{3i}}(\alpha_{3i-1},\alpha_{3i+1})\bigwedge_{1\leqslant j<\frac{m-1}3}\starr_{\beta_{3j+1}}(\beta_{3j},\beta_{3j+2}).
\end{equation}
Редукцию по этим же рёбрам производим и в склеенном графе $\Gamma'$ (это корректно).

Рассмотрим редукцию класса $v\circ\hat u(\hat\bigcirc)$, где $v$ обозначает расклейку, рассмотренную ранее, при которой происходит распад. Она совпадает с образом при $v$ редукции элемента $\hat u(\hat\bigcirc)$, который можно вычислять по формуле леммы \ref{raskleim-prod}, так как в выражении (\ref{raskleikacirc}) нетривиальной будет расклейка только одного множителя, а именно, $\hat u(\starr_{\hat\alpha}(\beta_{m-1}\alpha_1))=\starr_{\gamma}(\beta_{m-1}\alpha_1)$. Нетрудно убедиться, что редукция этого образа совпадает с редукцией элемента (\ref{netrivialniyv8}). Пусть (\ref{netrivialniyv8}) является расклейкой \S-неразложимого класса $[z']\in H'$. Тогда полученное утверждение означает, что $[z']-\hat u(\hat\bigcirc)\in\Ker v$, то есть $[z']_*=\hat u(\hat\bigcirc)_*$, так как $\Ker u=\Im v$ состоит из \S-разложимых элементов. Но \SS-класс элемента $[z']$ является порождающим алгебры $H'/(\S)$, значит, тоже самое можно сказать и о $\hat u(\hat\bigcirc)$.
}

\section*{Благодарности}
Выражаю благодарность Евгению Соломоновичу Голоду (1935--2018), под чьим непосредственным руководством была выполнена эта работа.

\bibliographystyle{unsrt}  
%\bibliography{references}  %%% Remove comment to use the external .bib file (using bibtex).
%%% and comment out the ``thebibliography'' section.

%%% Comment out this section when you \bibliography{references} is enabled.

\end{document}